\documentclass[preprint,12pt]{elsarticle}

\usepackage{graphicx}
\usepackage{amssymb}
\usepackage{amsthm}

\usepackage{lineno}

\usepackage{nameref}
\usepackage{array}
\usepackage{animate}
\usepackage{epstopdf}
\usepackage{verbatim} 
\usepackage{array} 
\usepackage{float}
\usepackage{mathpazo} 
\usepackage{changepage} 
\usepackage{amsfonts,amsmath,mathrsfs}
\usepackage{subcaption}
\usepackage{booktabs} 

\journal{Mathematical Biosciences}

\begin{document}

\begin{frontmatter}


\title{Effects of mutations and immunogenicity on outcomes of
  anti-cancer therapies for secondary lesions}



 \author[label1,label2,label4]{Elena Piretto } 
  \author[label1]{Marcello Delitala}
   \author[label3]{Peter S. Kim}
    \author[label4]{Federico Frascoli \corref{cor1}}
\cortext[cor1]{Corresponding author}
 \ead{ffrascoli@gmail.com}
 \address[label1]{Department of Mathematical Sciences, Politecnico di Torino, Turin, Italy}
 \address[label2]{Department of Mathematics, Universit\'a di Torino, Turin, Italy}
\address[label3]{School of Mathematics and Statistics, University of Sydney, Sydney, New South Wales, Australia}
\address[label4]{Department of Mathematics, Faculty of Science, Engineering
and Technology, Swinburne University of Technology, Hawthorn, Victoria, Australia}

\begin{abstract}
Cancer 
development is driven by mutations and selective forces, including the
action of the immune system and interspecific competition. When
administered to patients, 
anti-cancer therapies affect the development and
dynamics of tumours, possibly with various degrees of resistance due to
immunoediting and microenvironment. 
Tumours are able to
express a variety of competing phenotypes with different attributes
and thus respond differently to various anti-cancer
therapies. 

In this paper, a mathematical framework incorporating a system of delay
differential equations for the immune system activation cycle and an
agent-based approach for tumour-immune interaction is presented.
The focus is on those metastatic, secondary solid lesions that are
still undetected and non-vascularised. 

By using available experimental data, we analyse the effects of combination therapies
on these lesions and
investigate the role of mutations on the rates of success of common
treatments. Findings show that mutations, growth properties and
immunoediting influence therapies' outcomes in nonlinear and complex
ways, affecting cancer lesion morphologies, phenotypical compositions
and overall proliferation patterns. Cascade effects on final
outcomes for secondary lesions are also investigated, showing that actions on primary lesions could
sometimes result in unexpected clearances of secondary tumours. This
outcome is strongly dependent on the clonal composition of the primary
and secondary masses and is shown to allow, in some cases, the control
of the disease for years.  
\end{abstract}


\end{frontmatter}



\section{Introduction}
\label{Sec:intro}

Cancer is a generic definition of a disease that, among
its typical features, is driven by dynamic alterations in the
genome~\cite{hanahan2011hallmarks}.
These microscopic changes not only give
birth to a variety of different
types of cancer at the macroscopic scale, but can also lead to heterogeneity within the same
cancer tissue: tumour phenotypes undergo clonal expansion and
genetic diversification, promoting natural selection mechanisms that
favor cell clones
with advantageous characteristics~\cite{tabassum2015tumorigenesis, gerlinger2010darwinian}.

Alterations in the DNA of the cell, such as inclusions of copy number
aberrations and point mutations, occur early during the neoplastic
transformation and usually before any possible clinical
detection~\cite{sottoriva2015big}. The step-wise accumulation of
driver mutations may confer survival advantages in relation to the
particular environment in which they are embedded and may be accelerated by
so-called selective sweeps~\cite{rubben2017cancer}. Furthermore,
although the immune system
routinely recognises and kills any dangerous host including cancer,
mutations can provide cancer cells with the ability to avoid detection or
immuno-suppress the environment, advantaging tumour progression or
preventing eradication~\cite{hanahan2011hallmarks}.  
Processes involving mutations, cell
growth and immune surveillance cumulatively result in the emergence of different
cancer populations integrated in an environment made up of healthy tissue,
immune cells and stroma~\cite{pacheco2014ecology, hillen2014mathematical}.

Understanding how these complex interactions shape and influence each
other is one of the greatest challenges in current medical
biosciences. For example, morphology is known to be 
strongly sensitive to tumour adaptation to the environment (e.g. the lack of
nutrients, oxygen, space) and by the combined action of immune response 
and existing anti-cancer therapies~\cite{anderson2006tumor} such as
chemotherapy, radiotherapy, immune-boosting and so on. In the last quarter of century, a number of diverse
contributions have been proposed from the biomathematical community to
shed light on some of these complex interaction mechanisms.
 Several
mathematical models have been advanced using the framework of population dynamics,
with tumour immune interactions considered, for example, in Ref.~\cite{de2006mixed,
  wilson2012mathematical} and cancer mutations in
Ref.~\cite{asatryan2016evolution}. Other works have involved a discrete
Cellular Potts approach~\cite{sottoriva2011modeling} or different
degrees of hybrid
modelling~\cite{anderson2008integrative, jeon2010off}, with particular
focus on tumour
shape~\cite{anderson2006tumor,jiao2012diversity}. The effects of
some of the currently available anti-tumour therapies have also been analysed in the context of
evolutionary dynamics~\cite{gatenby2009adaptive, gatenby2009lessons},
with immunotherapy~\cite{eladdadi2014modeling, frascoli2014dynamical,
  piretto2018combination, bunimovich2008mathematical} and, recently, using agent-based modelling in the context of virotherapies~\cite{jenner2018modelling}. A number of reviews 
detailing the evolution and the contribution of these and
other models also exist in the
literature~\cite{bellomo2008, wilkie2013review,
  eftimie2011interactions, eladdadi2014mathematical}.

The focus of the present work is on metastatic secondary solid
lesions, with particular emphasis on the role of the immune system and
mutations. Scope of the this work is the study of the effects of
different combination therapies on secondary lesions in order to better
understand the dynamics involved and the role of mutations on 
treatments' effectiveness. The rest of the paper is organised as
follows. In the 
``\nameref{Sec:model}'' section, a description of the mathematical approach used to
describe tumours, immune responses and anti-cancer therapies is given. Findings obtained via
computational analysis are illustrated and analysed in the
``\nameref{sec:result}'' and  ``\nameref{sec:discuss}''  sections. Finally, the 
``\nameref{sec:conclusion}'' section terminates the paper.

\section{Model}
\label{Sec:model}
Let us consider the biological setting under study as follows: a
primary, clinically detected cancer is present in a patient and it is
scheduled to be treated with different therapeutic approaches, in an
effort to improve the patients' clinical outlook. A
secondary lesion is also growing, undetected and located away from the primary
site, due to previous metastatic events and migration of tumour
cells belonging to the first lesion. We are interested in
understanding how the secondary lesion is affected by strategies aimed
at reducing the primary one. {\it Our approach is based on an existing
mathematical model for tumour-immune interaction
\cite{kim2012modeling}, which has been validated previously both from
the point of view of biological appropriateness and sensitivity to
model parameters. The phenomena at hand are inherently complex and
there is a number of unknowns that still characterise these
processes. Our work is thus focussed on understanding the major trends
and the typical outcomes that can emerge in treating secondary
lesions, providing some quantitative data that can be tested
experimentally.}

The dynamics between a heterogeneous, small, solid
cancer lesion and the immune system is formulated using an hybrid agent-based model
(ABM) coupled with a delay differential equation (DDE) system. An 
immune response to cancer cells that grow and mutate is simulated using a population
of cytotoxic T lymphocytes (CTLs), which mature in a tumour-draining lymph node. The overall 
approach rests on an existing framework, originally
discussing tumour cells endowed with only a unique, single
phenotype. The novelty of the present formulation lies in considering 
more than one clone, with mutation processes strongly influencing and
shaping tumour growth dynamics.
For a full analysis and description
of the model we refer the reader to Ref.~\cite{kim2012modeling}, and
only discuss the equations briefly in the following.

\begin{figure}[htbp!]
\centering
\includegraphics[width=1\textwidth]{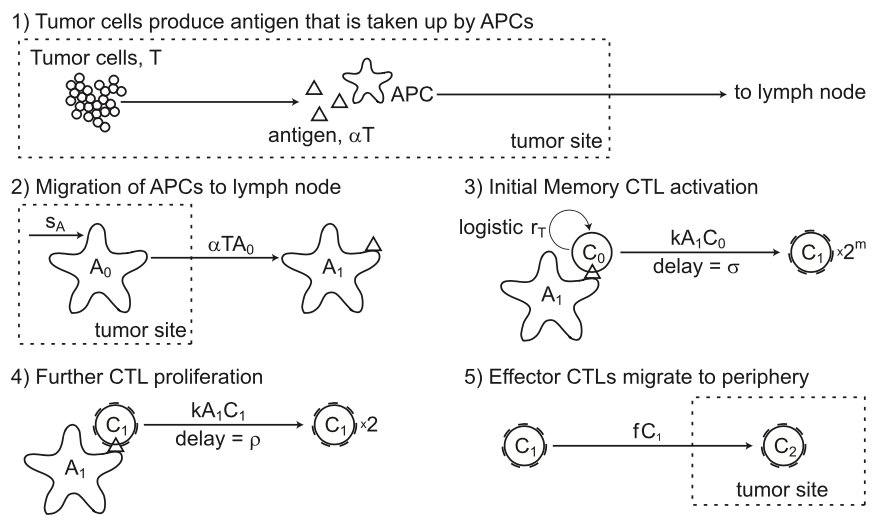}
\caption[Immune system activation cycle]{{\bf Immune system activation cycle} as described by the system of equation \eqref{eq:dde}. The behaviour of the immune system is modeled as ~\cite{kim2012modeling}}
\label{au:figimmuno} 
\end{figure}

The system describing immune activation is given by: 
\begin{equation}
    \begin{cases}
    A'_0(t) =& s_A-d_0A_o(t)-\alpha T(t)A_0(t), \\
    A'_1(t) =& V_{\mathrm{ratio}}\alpha T(t)A_0(t)-d_1A_1(t),\\
    C'_0(t) =& r_C \left( 1-\frac{C_0(t)}{K} \right) C_0(t)-\mu A_1(t) C_0(t), \\
    C'_1(t) =& 2^m \mu A_1(t-\sigma)C_0(t-\sigma)-\mu A_1(t)C_1(t) + \\
    & + 2 \mu A_1(t-\rho) C_1(t-\rho)-\delta_1 C_1(t)-f C_1(t),\\
    C'_2(t) =& \frac{f C_1(t)}{V_{\mathrm{ratio}}}-\delta_1C_2(t),
    \end{cases}
    \label{eq:dde}
\end{equation}
where $T$ is the total cancer cell population and $A_0$, $A_1$, $C_0$,
$C_1$, $C_2$ are the concentrations of antigen
presenting cells (APC), mature APCs, memory CTLs, effector CTLs and
CTLs, respectively. A sketch of the dynamics captured by the above
equations is depicted in Fig.\ref{au:figimmuno}.

The first two equations describe the transition from immature APCs
circulating in the periphery to mature ones migrating to the lymph
node as a response to tumour antigens. 

The population of immature
APCs is generated and dies at constant rates $s_A$ and
$d_0$, with the maximum value of $A_0$
corresponding to the equilibrium level $s_A/d_0$. When tumour antigens
are presented, $A_0$ decreases proportionally
to the antigenicity value $\alpha$ and 
mature APCs ($A_1$) begin entering the lymph node, with some dying at
natural death rate $d_1$. 

The presence of mature, tumour-antigen-bearing APCs in the lymph node
causes memory CTLs to activate and mature into effector CTLs, with a
certain delay. Consequently, the effector CTLs proliferate and
migrate to the tumour site where the anti-tumour immune response
starts. This process is captured as follows. The third equation represents the stimulation by the APCs of
the memory CTLs ($C_0$), with a logistic growth rate independent
of the external stimuli and a stimulation by
mature APCs that follows a mass action law. In the fourth equation, 
memory CTLs develop a minimal
division process, characterised by $m$ times divisions, and 
evolve in effector CTLs ($C_1$) with a time delay $\sigma$. Then, effector CTLs divide
again in a time $\rho$ and flow away of the lymph node or naturally
die with a rate $\delta_1$. The last equation represents the concentration of CTLs ($C_2$) in
the periphery around the tumour and provides the concentration $C_2$
used by the ABM component of the model to generate the
boundary conditions for the tumour-site domain. Table
\ref{tab:param_dde} shows parameter values used in this system of equations and their
meanings. Subsection \ref{sec:parsen} reports the parameter estimation and the related sensitivity analysis, which have been mostly performed in previous works ~\cite{kim2012modeling}.

\begin{table}[!ht]
\begin{adjustwidth}{-2cm}{-2cm} 
\begin{center}
\caption[Parameters used in the DDE component of the model]{{\bf Parameters used in the DDE component of the model. For more details about parameter estimation refer to ~\cite{kim2012modeling}}}\label{tab:param_dde}
\begin{tabular}[c]{|p{1cm}|p{4cm}|p{3.5cm}|p{0.75cm}|p{6cm}|}
\hline
{\bf Par.}  & {\bf Description } & {\bf Value (range) } &  {\bf Ref. }  &  {\bf Note } \\ \hline \hline

$A_0(0)$                & Initial concentration of immature APCs                    & 0.01 $\mathrm{k}/\mathrm{mm}^3$   & \cite{catron2004visualizing} & same order of magnitude as the APCs in the lymph node \\ \hline
$d_0$                   & Death/turnover rate of immature APCs                      & 0.03 $\mathrm{day}^{-1}$   & \cite{mohri2001increased} & similar to those of na\"ive T cells \\ \hline
$s_A$                   & Supply rate of immature APCs                              & $0.3 ~\mathrm{k}/\mathrm{mm}^3 \mathrm{day}^{-1}$   & \cite{mohri2001increased} & $s_A=A_0(0)d_0$ \\ \hline
$d_1$                   & Death/turnover rate of mature APCs                        & 0.8 $\mathrm{day}^{-1}$   & \cite{belz2007killer} & using a half-life of 20 h: $d_1=(\ln{2})/20 \; h^{-1}$  \\ \hline
$K$                     & Equilibrium concentration of memory CTLs          & $2\%\cdot 200 ~\mathrm{k}/\mathrm{mm}^3$  & \cite{catron2004visualizing}  & $2 \%$ of the T-cells in a lymph node of radius 1 mm  \\ \hline
$r$                     & Logistic growth rate of memory CTLs
                 & $\log 2~\mathrm{day}^{-1}$  & \cite{kim2012modeling} & minimum doubling time of 1 day \\ \hline
$m$                     & Minimal number of CTL divisions                           & 10    & \cite{wodarz2005effect} & range from 7 to 17 cell divisions \cite{de2001recruitment,kaech2001memory} \\ \hline
$\delta_1$              & Death/turnover rate of effector CTLS                      & 0.4 $\mathrm{day}^{-1}$    & \cite{de2003different} & half-life during T-cell contraction of 41 h: $\delta_1=(\ln{2})/41 \; h^{-1}$ \\ \hline
$\mu$                   & Mass-action coefficient                                   & 20 ($\mathrm{k}/\mathrm{mm}^3$)$^{-1}\mathrm{day}^{-1}$    & \cite{catron2004visualizing} &  $\mu = 0.5 \mu_0$ with $\mu_0=4.8 cell^{-1}day^{-1}$ and $V_{\mathrm{lymph \; node}}=8.4\cdot10^{-3}mm^3$\\ \hline
$\rho$                  & Duration of one CTL division                              & 8 h   & \cite{janeway1996immunobiology} & T-cell doubles every 8 hours during expansion \cite{de2003different} \\ \hline
$\sigma$                & Duration of CTL division program                          & $1+(m-1)\rho$   & \cite{veiga2000response} & first division does not occur until 24 hours after stimulation\\ \hline
$\alpha$                & Antigenicity of the tumour                                 & $10^{-9}$ ($\mathrm{k}/\mathrm{mm}^3$)$^{-1}\mathrm{day}^{-1}$  & \cite{kim2012modeling} &  reciprocal of the rate of encountering of antigen from a tumor cell by APC \\ \hline 
$f$                     & CTLs flow rate out of lymph node to tissue & $0.7 \; \mathrm{day}^{-1}$ & \cite{kim2012modeling} &  effector CTLs emigrate at a half-life of 1 day: $f=\ln{2} day^{-1}$ \\ \hline
$V_{\mathrm{ratio}}$    & Ratio of volume of tissue to the lymph node               & 1000   &  \cite{ying2009three} & lymph node compartment is $\sim 1 ml$ and the breast tissue $\sim 1 l$  \\ 
\hline
\end{tabular}
\end{center}
\end{adjustwidth}
\end{table}

\subsubsection*{Cancer dynamics and immune cells response}\label{subsec:ABM}
The ABM controls tumour growth dynamics and the
interaction between tumour cells and CTLs at the tumour site, which follow
specific algorithmic rules. Our model does not consider healthy
tissue around tumours and other structures such as the stroma or the
cells part of a vascular network: it is assumed that tumour's
surroundings are effectively healthy cells being ``pushed away'' by the
growing tumour. Note also that no vascularisation is present due to the limited
size of the secondary tumour lesion, which is considered to be small, solid and with no
necrotic core. Furthermore, other motility of metastatic processes from the secondary lesion are neglected. The overall assumption is that the secondary tumour is trying to colonize the site and is in its early stages of proliferation. 
All cells partaking the dynamics are represented as spheres of radius $r$ in 3D space, with no overlap.
ABM is updated in discrete timesteps $\Delta t$.

{\bf CTL agents}. The rules that govern CTLs cells via the ABM are
three: motion around the tumour, recruitment of other immune system
cells and killing of tumour cells.  As mentioned, CTLs cells appear at a
concentration $C_2(t)$ at the border of
the spherical domain representing the region of interest where the tumour is
growing. They then move into that region performing Brownian motion in
3D space until they either collide with a cell or leave the
domain. At each time step, the position of the
cells are given by independent random variables with normal distribution
$\mathcal{N}(0,\sigma^2 \Delta t)$, where the variance is such that
$\sigma^2=2D$, with $D$ being the diffusion rate of the CTLs. When an
immune cell comes into contact with a cancer cells three possibilities
exist: 
\begin{itemize}
    \item A CTL clone can be recruited with a probability $1-e^{-\Delta
        t/C_{\mathrm{recruit}}}$, with $C_{\mathrm{recruit}}$ being
      the average recruitment time. Mathematically, CTL recruitment is modeled similarly to Mallet et al. \cite{mallet2006cellular} with cellular automata, and it is biologically validated as in \cite{soiffer2003vaccination, soiffer1998vaccination}. When the first CTL cell
      engages a cancer cell and starts recruiting another CTL clone, a
      second cell appears at a position adjacent to the first
      cell. The direction of the new clone is chosen randomly among
      all directions available. 

   \item A cancer cell is not recognised with a probability $1 -
      P_{i,\mathrm{recog}}\cdot\Delta t$, where $P_{i,\mathrm{recog}}$ is
      the probability of the i-th cancer phenotype (see below) to be
      recognised by the immune system. The parameter
      $P_{i,\mathrm{recog}}$ has a value of one for cancer agents, expressing antigens completely matching with the T-cell receptors and thus, that are always
      detected by the immune system. A value of zero indicates that
      the antigens of a phenotype are completely unrecognised. If the
      cancer cell is not
      recognised, the CTL starts to move again choosing a new random direction and 
      accelerating up to the maximum
      unit standard deviation $\sigma_{\mathrm{max}}$. If
      $C_{\mathrm{acc}}$ is the time necessary to accelerate from the
      stationary to the maximum diffusion rate, the CTL acceleration
      is computed as: $\sigma(t)=\sigma_{\mathrm{max}}\cdot
      \mathrm{min} \left( t/C_{\mathrm{acc}},1 \right)$. 
 This approach aims at approximating CTL chemotaxis along a chemokine 
gradient \cite{mackay1996chemokine, maurer2004macrophage}.

    \item A cancer cell is recognised and killed with a probability
      $1-e^{-\Delta t/C_{\mathrm{kill}}}$, with
      $C_{\mathrm{kill}}$ being the average time for a CTL to eliminate
      a cancer cell. The killing process is obtained
      by removing the agent. After the agent is removed, the
      immune cell starts to move again as described above.
\end{itemize}
If CTLs die naturally, then they are removed from the system. An
explanation of the ABM-parameters is reported in Table
\ref{tab:param_abm} whereas parameter estimation and sensitivity analysis is discussed in Subsection \ref{sec:parsen}.

\begin{table}[!ht]
\begin{adjustwidth}{-.6in}{-.6in} 
\begin{center}
\caption[Parameters used in the ABM component of the model]{{\bf Parameters used in the ABM component of the model.  For more details about parameter estimation refer to ~\cite{kim2012modeling}}}\label{tab:param_abm}
\begin{tabular}[c]{|p{1.1cm}|p{5cm}|p{2cm}|p{0.75cm}|p{6.5cm}|}
\hline
{\bf Par.}  & {\bf Description } & {\bf Value (range) } &  {\bf Ref. }  &  {\bf Note } \\ \hline \hline

$\Delta t$              & Time step                                                 & 1 min   & \cite{kim2012modeling}  & timescale of the fastest dynamic simulated in the model  \\ \hline
$r$                     & Radius of cells                                           & 5 $\mu \mathrm{m}$   & \cite{lin2004model}  &  \cite{mallet2006cellular, alarcon2003cellular, catron2004visualizing} \\ \hline
$T_{\mathrm{div},i}$    & Avg. division time of i-th cancer phenotype               & 1-39 day    & \cite{kirschner1998modeling}  &  \cite{mallet2006cellular, kuroishi1990tumor, michaelson2003estimates} \\ \hline
$T_{o}$                 & Avg. division time of {\it original} phenotype            & 7 day    & \cite{kirschner1998modeling}  &  \cite{mallet2006cellular, kuroishi1990tumor, michaelson2003estimates}  \\ \hline
$\sigma_{\mathrm{max}}$ & Max unit standard deviation of CTL diffusion              & 12 $\mu \mathrm{m} ~\mathrm{min}^{-1}$   &  \cite{catron2004visualizing} & \cite{friedl2001interaction}  \\ \hline
$C_{\mathrm{acc}}$      & CTL acceleration time from 0 to $\sigma_{\mathrm{max}}$   & 5 h     & \cite{kim2012modeling}  &   \\ \hline
$C_{\mathrm{death}}$    & Avg. CTL lifespan                                         & 41 h    & \cite{de2003different}  &   \\ \hline
$C_{\mathrm{recruit}}$  & Avg. time fro CTL recruitment                             & 22 h   &  \cite{kim2012modeling} &  \cite{soiffer2003vaccination, soiffer1998vaccination} \\ \hline
$C_{\mathrm{kill}}$     & Avg. time fro CTL to kill tumour cell                      & 24 h    &  \cite{kim2012modeling} & killing target cells may require a long recovery period   \\ \hline
$R$                     & Radius of region of interest                              & 620.4 $\mu \mathrm{m}$   &  \cite{kim2012modeling} &   \\ \hline
$h$                     & Thickness of CTL cloud                                    & $3\sigma_{\mathrm{max}}\sqrt{\Delta t}$    & \cite{kim2012modeling}  & probability that a CTL could pass from outside into the region of interest is 0.001  \\ \hline  
$P_{\mathrm{mut}}$      & Probability of mutation                                   & 0.01 $\mathrm{min}^{-1/2}$   & \cite{wood2007genomic}  &   \\ \hline
$P_{\mathrm{recog},i}$  & Probability of recognition of i-th cancer phenotype       & 0-1    &  &  Span the entire probability range \\ \hline
$P_{o}$  & Probability of recognition of {\it original} phenotype       & 1   &   &  the APC cell can always recognize the antigen released \\ 
\hline
\end{tabular}
\end{center}
\end{adjustwidth}
\end{table}

{\bf Cancer agents}. Tumour cells can
proliferate, mutate or die, killed by the immune system, and no migration
is considered. This approximation is motivated by the scope of the
study, which is focused on the solid, growing secondary lesion after
the colonisation of a new tissue. In this early stage of implantation
most of the cells are assumed to be in a proliferation state and
migration can be neglected \cite{hatzikirou2012go}.
Cellular division occurs with a probability $1-e^{-\Delta t/T_{i,
    \mathrm{div}}}$, $i=1,...,5$, where $T_{i, \mathrm{div}}$ is the
average division time of the i-th tumour phenotype. When a tumour cell
divides, the position of a new cell is chosen randomly on the mother
cell's perimeter, such that the daughter cell is tangent. If no space
is available in the chosen position, the division process
fails and no new agent is created, mimicking the contact-inhibition
mechanism occurring in the early stages of metastasis implantation \cite{mendonsa2018cadherin}.    

To analyse the effect of mutations on cancer development and immune
response, we use five different cancer phenotypes that may emerge from
the mutation of an {\it original} clone, identified by different values of
characteristic parameters $T_{\mathrm{div},o}=T_o$ and $P_{\mathrm{recog},o}=P_o=1$. 
Mutations can occur during cell duplication, 
with a probability $P_{\mathrm{mut}}\cdot \Delta t$ that aims to
capture the genetic instability of the system. 
Each mutated cell is then identified by indices
representing the level of expression of the two characteristic quantities
$T_{\mathrm{div}}$ and $P_{\mathrm{recog}}$. These values effectively classify the
mutated clones and the following mutated phenotypes. 
Modeling few phenotypes of mutated cells is a simplification justified
by several works showing that only a limited number of phenotypes are
predominant in a tumour, see for example \cite{anderson2006tumor}. For the
scope of our study, the five mutated clones are prototypical of a wide
range of similar mutations. 
In Table \ref{tab:clones} cancer clonal composition is considered. One of the assumptions is that only one class of CTLs is modeled and it is not antigen specific. Although different types of CTLs could take part in an immune response and act differently depending on the clone, our immune attacks are regulated only via $P_{recog,i}$.

\begin{table}[!ht]
\begin{adjustwidth}{-.6in}{-.6in} 
\begin{center}
\caption[Cancer clonal composition]{{\bf Cancer clonal composition}}\label{tab:clones}
\begin{tabular}{|c|l|c|c|}
\hline
{\bf Name}  & {\bf Description } & {\bf $P_{\mathrm{recog}}$}  & {\bf $T_{\mathrm{div}}$ } \\ \hline \hline
{\it original}   & \begin{tabular}{@{}l@{}}
              First metastatic breast cancer clone  \\
              that colonises the new tissue, with \\
              evolutionary potential of phenotypic mutations. \end{tabular} 
              & 1 $\mathrm{min}^{-1}$  & 7 day \\
              \hline
(0.5,0.5)   & \begin{tabular}{@{}l@{}}
              Clone proliferates at 
              the same rate of the {\it original} \\ clone, 
              but has an increased ability to hide from \\ 
              the immune system.   \end{tabular} 
              & 0.5 $\mathrm{min}^{-1}$  & 7 day \\
              \hline
(0,1)       & \begin{tabular}{@{}l@{}}
              Clone is not recognised by CTLs, but \\
              the evolutionary cost of its ability is\\
              paid in term of proliferation: this \\
              phenotype is the slowest to reproduce.      \end{tabular} 
                & 0 $\mathrm{min}^{-1}$ & 13 day\\
              \hline
(0.25,0.75) & \begin{tabular}{@{}l@{}}
              Clone has intermediate properties: \\
              strong ability to hide and slow  proliferation   \end{tabular} 
               & 0.25 $\mathrm{min}^{-1}$ & 10 day \\
              \hline
(0.75,0.25) & \begin{tabular}{@{}l@{}}
              Clone has intermediate properties: \\ 
              weak ability to hide and fast proliferation   \end{tabular} 
              &  0.75 $\mathrm{min}^{-1}$ & 4 day  \\
              \hline
(1,0)       & \begin{tabular}{@{}l@{}}
              Clone has the ability to reach \\
              high number of cellular duplication, but is  \\
              always recognised by the immune system.   \end{tabular} 
               & 1 $\mathrm{min}^{-1}$ & 1 day \\
              \hline
\end{tabular}
\end{center}
\end{adjustwidth}
\end{table}

Using these different types of phenotypes, as we will see shortly,
helps us to shed light on the role of mutations in determining the
effectiveness of immune response and anti-cancer therapies. 
Different clonal compositions and reproductive and immunoediting
advantages dramatically influence the outcomes of anti-cancer therapies.

\subsection*{Modeling therapies: chemotherapy, immune boosting and radiotherapy} \label{sec:fitland}

One of the typical features of secondary lesions is that they usually
show cells with mutated
functional characteristics respect to the original tumour, due to the genetic instability typical of
metastatic masses they originate from. 
We reiterate that there is no analysis of the fate
of the global cancer disease but only on such secondary lesions, which
can show different dimensions, compositions, structures and biological
characteristics from the primary neoplasia. Chemotherapy, immune
boosting and radiotherapy are the strategies our modelling focuses on.

{\bf Chemotherapy.} This treatment consists of cytotoxic drugs
targeting a specific cellular
phase of the cell cycle to induce cell death. The procedure acts
against rapidly proliferating cells, independently from their nature 
~\cite{brunton2011therapeutics}. This means that healthy cells and
immune system cells are usually damaged along with cancer cells, and
this leads to well-known side effects for the patients. In this work,
only the primary killing effect against cancer cells and no
direct effects on the immune system is assumed. This simplification is
motivated by two main points. First, the average CTL lifespan is 41
hours, whereas the tumour division rate is greater and the tumour death
rate due to the therapy is slower. CTL cells are rapidly affected by
the reduction of the tumour mass and no new CTL is recruited: the
``old" cells tend to naturally die. Second, if on one side
chemotherapy affects the immune cells, on the other specific T-cell
response is reinforced \cite{zitvogel2008immunological}, and 
the investigation of these secondary effects is not in the scope of
the future present work.

During a cycle of chemotherapy of duration $Ch_{\mathrm{time}}$, the
i-th cancer phenotype can go through cellular death with probability
$1-e^{-{\Delta t}/{Ch_{\mathrm{kill},i}}}$, where the average time for
the drugs to induce cellular death is $Ch_{\mathrm{kill},i}$ and
depends from the proliferation potential of the phenotype. $Ch_{\mathrm{time}}$ takes into account a single cycle of three injections every three days and represents the global time duration of the chemotherapy's effects. The drug remains two days above a certain percentage level such that the cytotoxic effects on tumour cells can be considered constant. 

Different values of $Ch_{\mathrm{time}}$ have been explored as reported in Table \ref{tab:paramTh}, supposing that the same total dose is inoculated in continuous cycles of low metronomic doses. The effect of different $Ch_{\mathrm{time}}$ with the same total dose is a faster or slower decrease of the cancer population with similar qualitative dynamics.
In particular, for clone $(0,1)$ (refer to Table \ref{tab:clones} for notation), i.e. the phenotype that grows slowly
but is poorly immunogenic,
$Ch_{\mathrm{kill},\mathrm{(0,1)}}=Ch_{\mathrm{time}}$, namely a
$(0,1)$-death is very rare. Clone $(0+0.25i,1-0.25i)$, with
$i=1,...,4$, has
$Ch_{\mathrm{kill},\mathrm{(0+0.25i,1-0.25i)}}=Ch_{\mathrm{time}}-i\cdot
Ch_{\mathrm{eff}}\cdot Ch_{\mathrm{time}}$, so that the tumour with
higher proliferation rate has very high probability to die due to the
effect of the drug.

{\bf Immune boosting.}
We use this generic term to capture the number of clinically available strategies that
potentiate an immune response. For example, a
treatment that is increasingly used for cancer patients is the so-called
adoptive cell transfer (ACT), where patients' own immune cells are
stimulated and modified to treat their tumour. There are several types of
ACTs that go under different acronyms depending on the boosting
strategy employed, with the most used ones nowadays being
TIL (Tumour infiltrating lymphocytes), TCR (Tumour cell receptors)
T-cell and CAR (Chimeric antigen receptors) T-cells
treatments~\cite{rosenberg2008adoptive}.

We concentrate in particular on TIL therapy, where T-cells are
extracted from the patient's tumour, grown {\it in vitro} to boost
their numbers and injected back into the patient to contrast cancer
progression. This strategy appears to be, for example, one of the most effective treatment
against metastatic melanoma~\cite{rosenberg2008adoptive}. In our
approach, TIL is modeled as a continuous increase of the CTLs concentration in
the cloud, depending on the value of $C_2$ at the starting time for
the therapy. The net increase is modeled by a number of $Bo_{\mathrm{eff}}$
cells for a short time $Bo_{\mathrm{time}}$. For simplicity, in
the following we refer to this treatment as immune boosting or
simply boost. The parameters used to model boost and
chemotherapy are explained and collated in Table~\ref{tab:paramTh}.
\begin{table}[!ht]
\begin{adjustwidth}{-.6in}{-.6in} 
\begin{center}
\caption[Parameters used to model therapies]{{\bf Parameters used to model therapies}. Ranges indicate that different therapies (single and combined) are simulated with different values.}\label{tab:paramTh}
\begin{tabular}{|c|l|c|}
\hline 
{\bf Parameter}  & {\bf Description } & {\bf Value (range) } \\ \hline \hline
$Ch_{\mathrm{time}}$    & Duration of a chemotherapy cycle  & 10-50 day    \\
$Ch_{\mathrm{eff}}$     & Effect of chemotherapy            & $0-1/4$        \\
$Bo_{\mathrm{time}}$    & Persistence time of boosting (TIL)     & 3 day  \\
$Bo_{\mathrm{eff}}$     & Number of CTL cells injected      & 500-1000     \\
\hline
\end{tabular}
\end{center}
\end{adjustwidth}
\end{table}

{\bf Radiotherapy.} Radiotherapy (RT) uses ionising radiation to
induce cell death in a localised area under treatment. This therapy has
several positive and negative feedbacks on the immune system, modulating
different compartments of the tumour microenvironment. In particular,
tumour-specific antigens and immune-stimulatory signals are released by
the dying cancer cells. 

Because of its contributing
primarily to the original, metastatic neoplasia, the effect of RT is
here modelled as an indirect effect
on the secondary lesion and is accounted for as 
as a restoring factor in the ability of CTL cells to recognise and kill various cancer phenotypes. 

\subsection{Parameter estimation and sensitivity analysis} \label{sec:parsen}

The biological significance of parameters and processes that underpin
the present model has been discussed at length elsewhere
\cite{kim2012modeling, Frascoli2017215}. In some cases, such as, for
example, parameters used for tumour
division time or cell radius, well-established values in the literature
have been used \cite{lin2004model, alarcon2003cellular,
  kirschner1998modeling, mallet2006cellular}. In other cases,
estimations from the available experimental and theoretical data have
been carried out.

A sensitivity analysis has originally also been carried out 
for eight parameters of the model: $T_{\rm div},~\sigma_{\rm
  max},~C_{\rm acc},~C_{\rm recruit},~C_{\rm kill},~K,~\mu, ~m$ and
$\alpha$. Other parameters have not been considered because their role is
known to be marginal. For instance, the replenishment rate for memory
CTLs is known to be irrelevant, since only a very tiny fraction of
memory CTLs ($~ 1\%$) is known to be affected by the tumour. Similarly,
the duration of CTL division (time delay parameter $\rho$ in the DDE)
is too small to impact the CTL division program as a whole and does
not influence final outcomes. Using Spearman's rank-order
correlations, tumour populations' values and extinction times, Kim et al.
have concluded \cite{kim2012modeling} that tumour division times $T_{\rm div}$, antigenicity
$\alpha$ and the number of divisions of memory CTLs upon activation
$m$ are the most sensitive parameters.

In this work, we use the same parameters proposed in the original paper, with the
only difference of $T_{\rm div}$, still chosen in the proposed
interval but capturing a more aggressive tumour (i.e. $T_{\rm
  div} = 1-39$ days). The effect on simulations is to shorten the
proliferating phase, which occurs at a larger growth rate and allows
for a quicker immune response. 

The new probability coefficients introduced here, i.e. $P_{\rm mut}$ and
$P_{\rm recog,~i}$, have different effects. By using $5$ different
simulations with different initial random seeds and $10$ different
values of the parameters, we conclude that $P_{\rm mut}$ has no effect
on the final outcomes of the system, but only accelerates or delays
the identical dynamics shown by the model. 

$P_{\rm recog,~i}$ instead has a notable effect on the system. When
$P_{\rm recog,~i} = 0$ for a given $i$-th clone, the immune system is
unable to eradicate that particular phenotype and, if no external
therapy is present, the tumour endlessly grows.

As far as the values for therapies' parameters are concerned,
i.e. $Ch_{\rm time}$, $~Ch_{\rm eff}$, $~Bo_{\rm time}$, $~Bo_{\rm
  eff}$, they are chosen so that the dynamics
between tumour, immune system and therapies display interesting
behaviours and does not result in an immediate negative or positive
outcome. 
In particular, $Ch_{\rm eff}$ and $Bo_{\rm time}$ have been
varied in a number of different instantiations of the model, with only
the cases $Ch_{\rm eff} = 0.25$, $Bo_{\rm time} = 1000$ cells used in
the discussion of results. Variations of those parameters do not alter
in a significant way the prototypical dynamics that we will discuss
shortly. Note that $Ch_{\rm eff} = 0.25$ has been chosen so that the
cytotoxic drug targets fast proliferating cells.

\subsection{Morphological and complexity measures} \label{sec:measures}

Three indices that capture the shape and cellular compositions of the
tumour mass are
introduced and monitored in our computational experiments. {\it Note that
these indices can guide the evaluation of collective
properties of the evolving tumours. They are useful to discriminate
between different evolutions of the cancer masses and have also been
validated in some in vitro experiments, as shown by other authors in
previous works \cite{sottoriva2010cancer, jeon2010off}.}

{\bf Roughness.} Although random proliferation of a group of
cells leads to an almost smooth and spherical object, a
tumoural mass with diverse clonal families under the action of the
immune system can present itself as a rough
aggregate. To account for this, a measure of
roughness $M$ is introduced, as the ratio
between the surface $S$ and the volume $V$ of the aggregate~\cite{sottoriva2010cancer}.
The minimum ratio is represented by a sphere $S_s/V_s=(4\pi
R_s^2)/(4/3\pi R_s^3)=3/R_s$, where the value has been
non-dimensionalised as follows:
$M_{\mathrm{min}}=\sqrt{S_s}/\sqrt[3]{V_s}=\sqrt{4\pi}/\sqrt[3]{(4/3)\pi}$. The
roughness index $M$, expressed in terms of the minimal ratio for a
sphere, is given by: 
\begin{equation}
    M=\frac{\sqrt{S}}{\sqrt[3]{V}}\cdot\frac{1}{M_{\mathrm{min}}}=\frac{\sqrt{4\pi S}}{\sqrt[3]{3(4\pi)^2V}}.
    \label{eq:rough}
\end{equation}
A compact, non-infiltrated, almost spherical tumour mass has an index
$M$ close to unity while a tumour with highly irregular borders, for instance a solid
tumour with fingers and clusters of invasive cells or a mass highly
infiltrated by the immune system, displays a higher value. 

{\bf Radius of gyration.} This value represents the radius of a
sphere that contains the whole tumour aggregate and reads:
\begin{equation}
R_g=\sqrt{\frac{\displaystyle\sum_{i=1}^{N_c} (\mathbf{r}_i-\mathbf{r}_{cm})^2 }{N_c}},
\label{eq:Rg}
\end{equation}
where $N_c$ is the total number of cancer cells and $\mathbf{r}_i$ is
the distance of each clone from the center-of-mass of the tumour
($\mathbf{r}_{cm}$) that can vary during the tumour progression.

{\bf Shannon Index.} This indicator is introduced to account
for the presence of different phenotypes within a tumour, with regards
to tumour heterogeneity and 
relative frequency of each clonal family ($p_i$). The Shannon index
$H$ is thus defined as:
\begin{equation}
H=-\frac{\displaystyle\sum_{i=1}^{s}p_i \ln(p_i)}{\ln(s)}
\label{eq:shannon}
\end{equation}
where $p_i$ is the relative abundance of the phenotype $i$ and $s$ is
the total number of different phenotypes (in our case $s=6$). For
simplicity, $H$ is then normalised to the interval $[0,1],$ where zero
indicates a homogeneous population with only one clonal family and unity
represents a fully heterogeneous population where all phenotypes
are equally present.
\\

\section{Results}\label{sec:result}
The model outlined in the previous sections is the basis for {\it in
  silico} experiments, where a different number of therapies and their
combinations are tried out for significant values of the parameter
set. Depending on the initial conditions, the system exhibits
three typical behaviours, namely eradication, sustained (irregular) oscillations
or exponential, uncontrolled growth when an immune response to a
growing tumour is present. We consider parameter values where the effect of clinical therapies are relevant. Cases where the tumour grows too fast or too slow, making the effects of therapies not noticeable, are excluded from our analysis. Stationary behaviour has
never been observed. Outcomes also
depend on the characteristics of cellular phenotypes present in the
growing mass, strongly influencing its speed of growth, its ability to
counteract the action of T-cells with immunoediting and its
morphological qualities, which can hinder the ability of the immune
system to effectively erode the cancer. 

\begin{figure}[htbp!]
\centering\noindent
\makebox[\textwidth]{\begin{subfigure}{0.42\textwidth}
\includegraphics[width=1\textwidth]{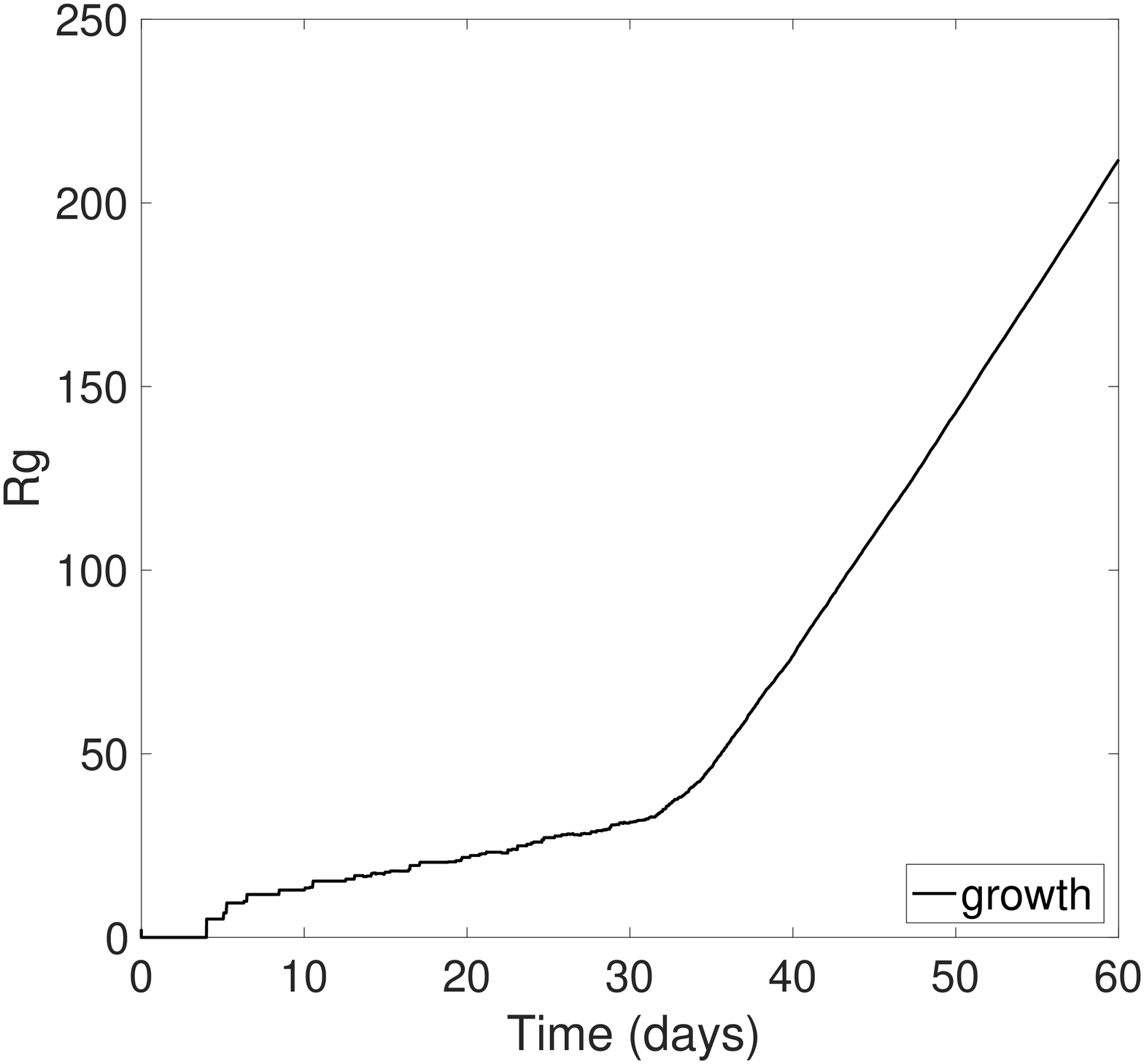}
	\caption{}
\end{subfigure}
\begin{subfigure}{0.42\textwidth}
\includegraphics[width=1\textwidth]{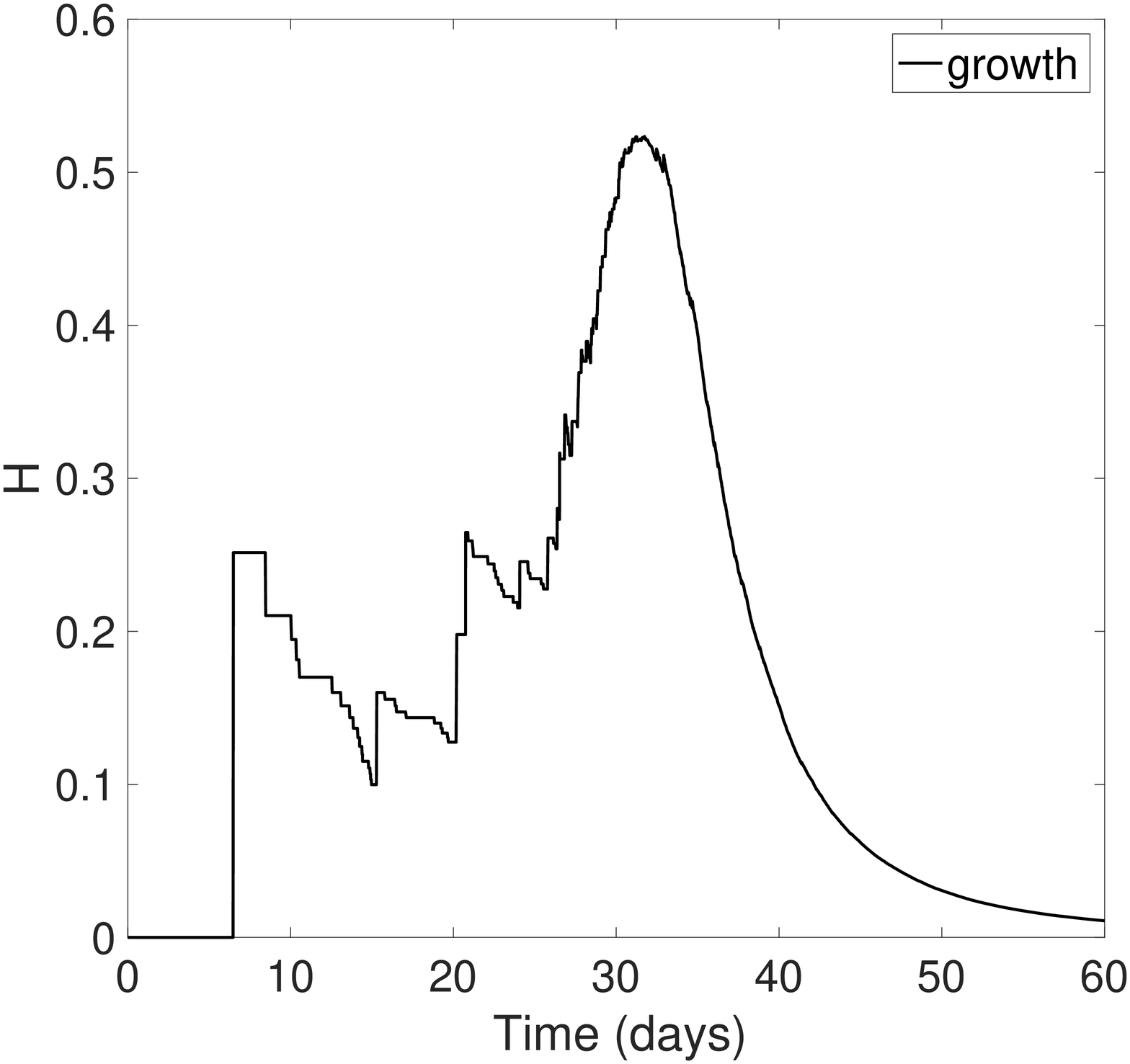}
	\caption{}
\end{subfigure}
\begin{subfigure}{0.42\textwidth}
\includegraphics[width=1\textwidth]{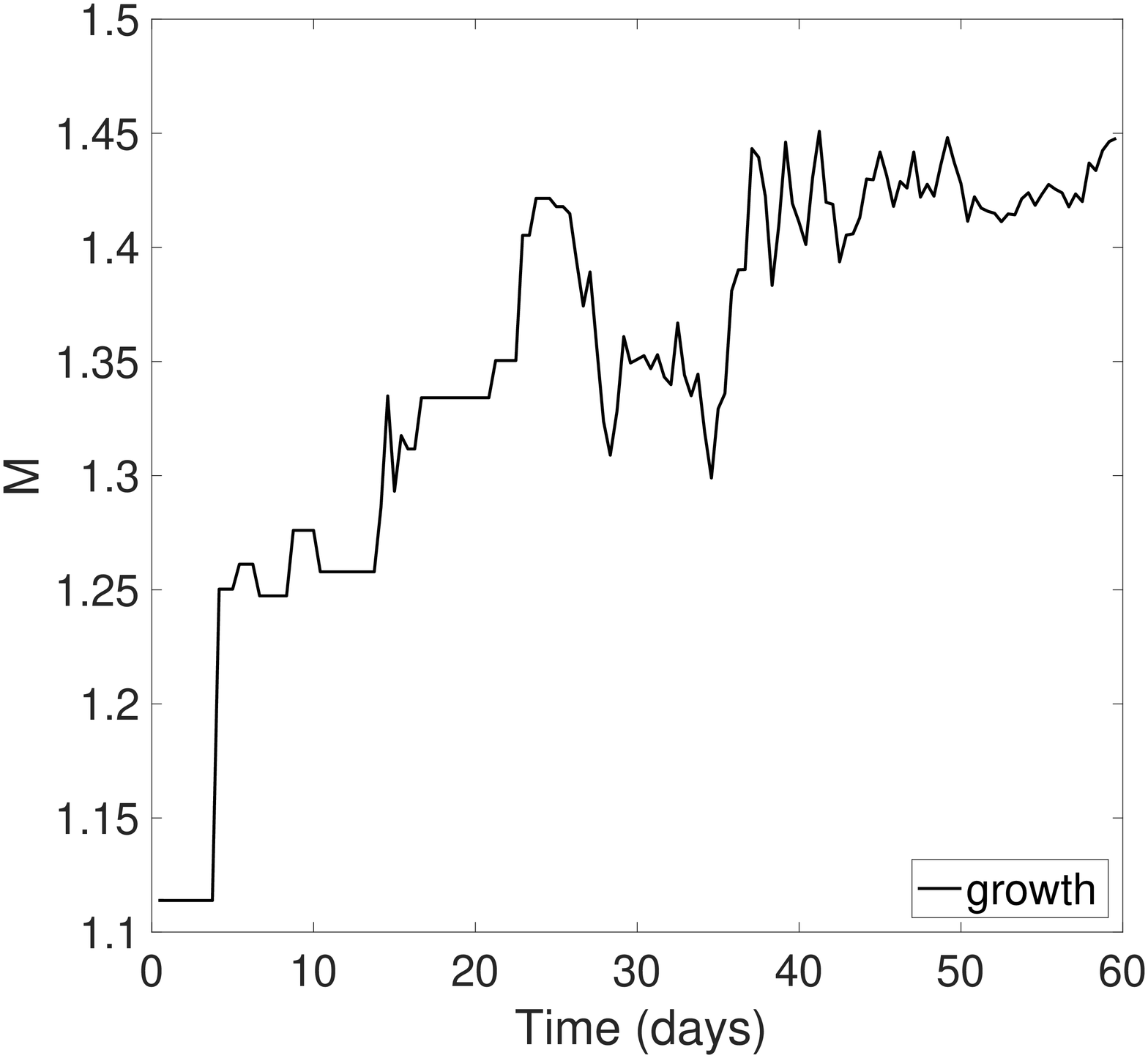}
	\caption{}
\end{subfigure}}\\

\makebox[\textwidth]{\begin{subfigure}{0.42\textwidth}
\includegraphics[width=1\textwidth]{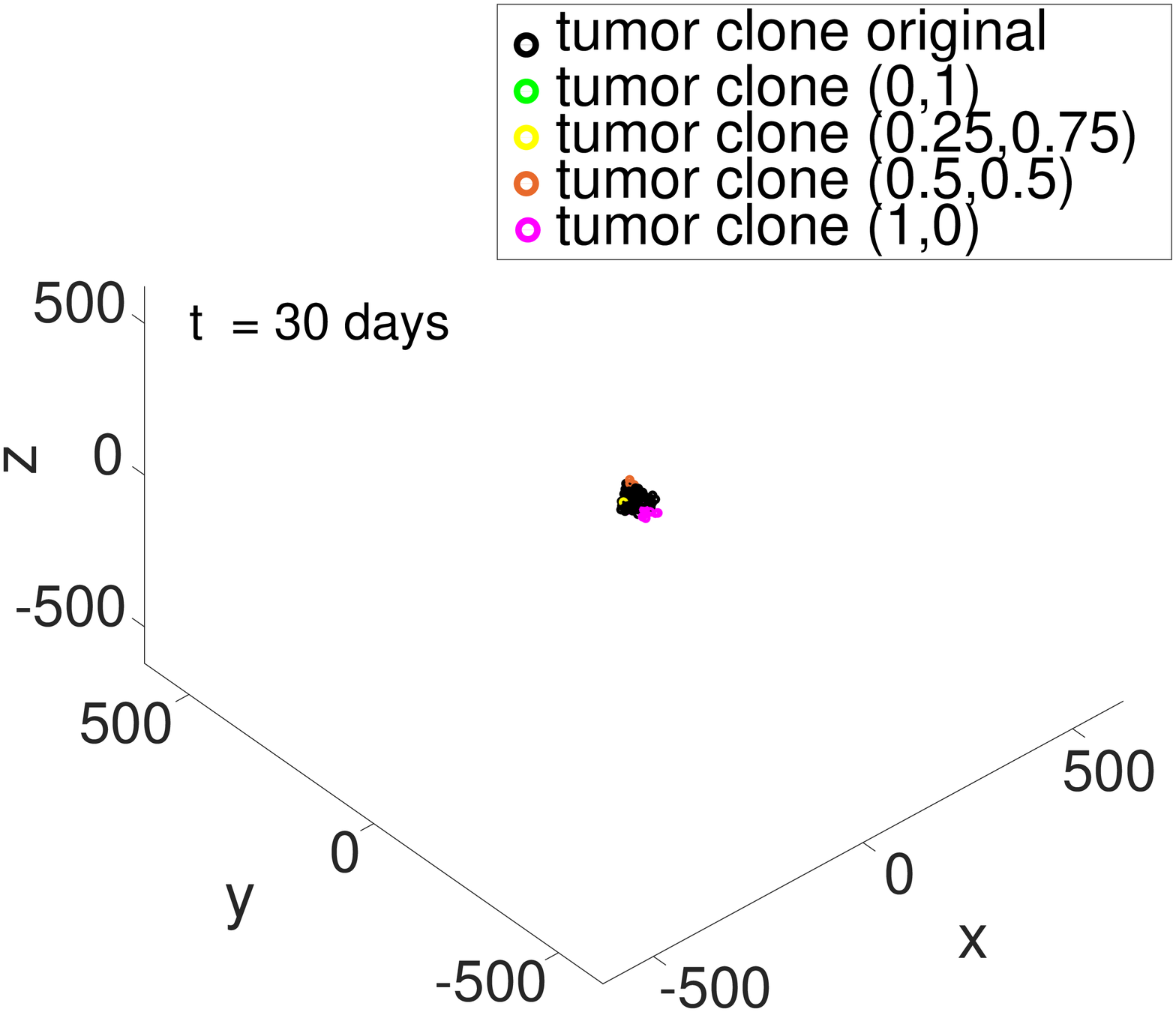}
	\caption{}
	\end{subfigure} 
\begin{subfigure}{0.42\textwidth}
\includegraphics[width=1\textwidth]{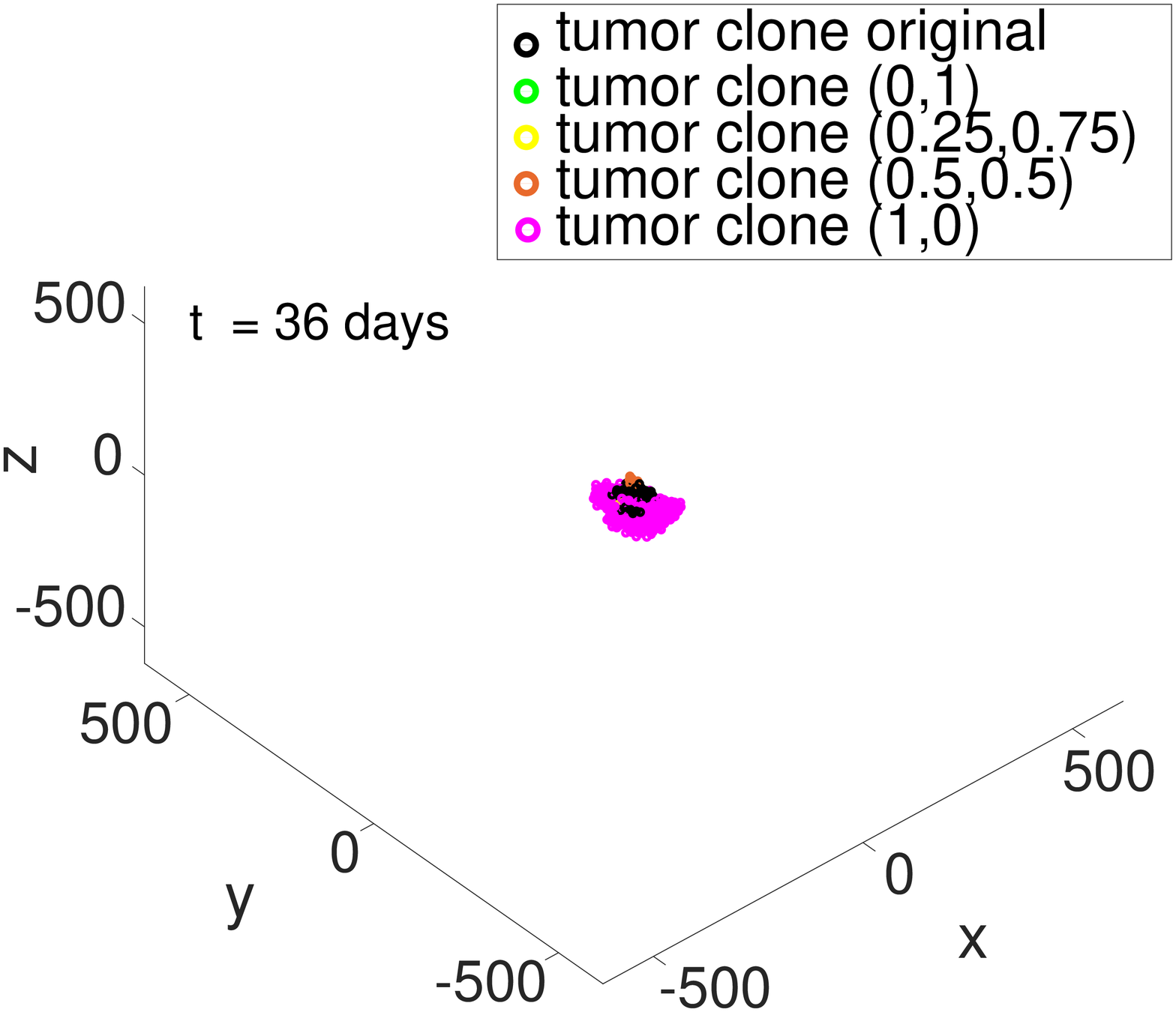}
	\caption{}
\end{subfigure}
\begin{subfigure}{0.42\textwidth}
\includegraphics[width=1\textwidth]{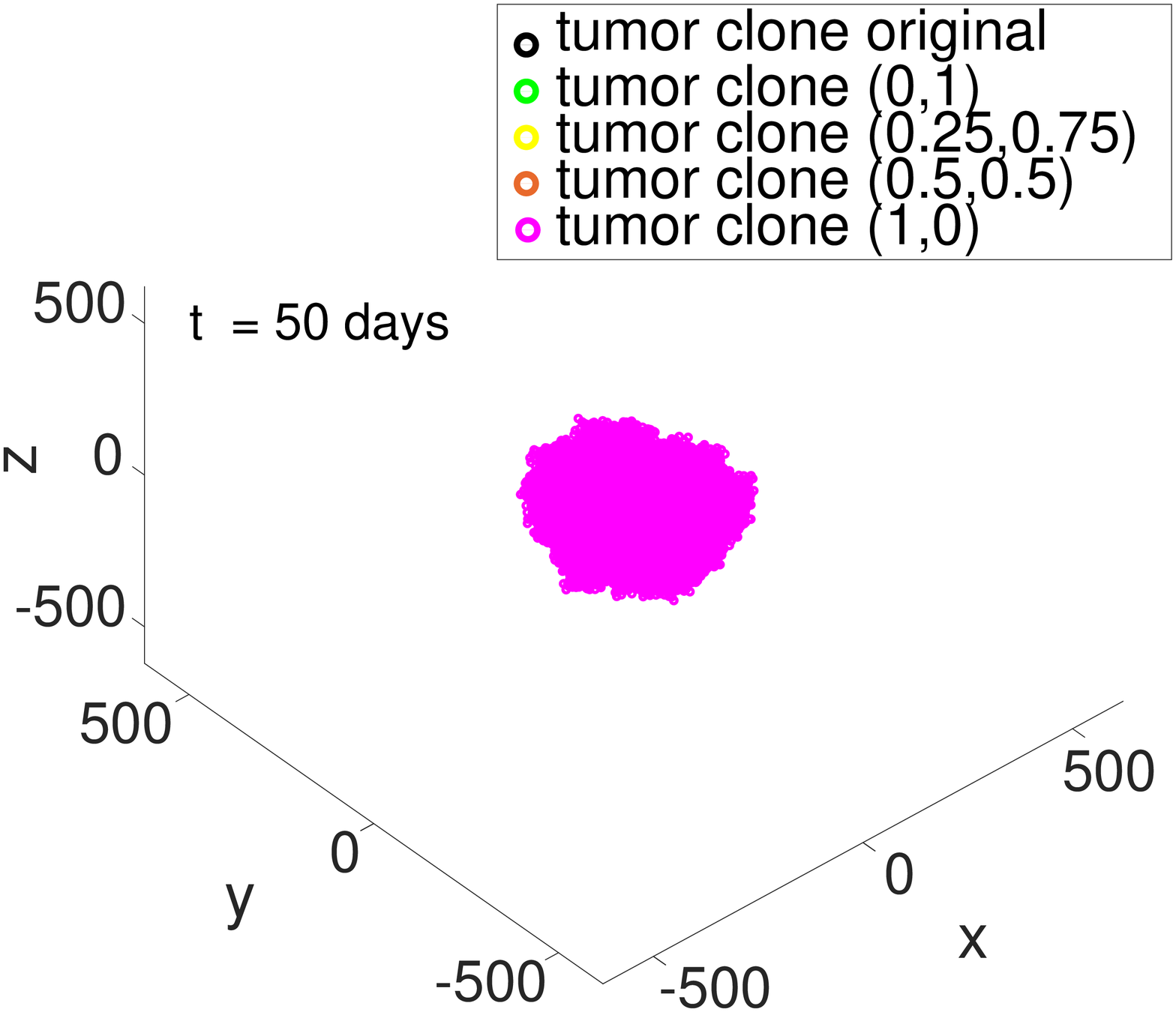}
	\caption{}
\end{subfigure}}

\caption[Tumor growth dynamic]{{\bf Tumour growth dynamics}
Panel (A) $R_g$ as function of time, slow initial trend followed by a
fast linear growth.  Panel (B) $H$ as function of time, increases with
mutations and decreases when the fittest clonal population outnumbers
the others. 
Panel (C) $M$ as function of time. The roughness index increases with
mutations, 
is higher when a new population with fast proliferation arises and 
reaches an almost stationary level when $H$ is near zero.
Panel (D), (E), (F) 3D-views of the tumour, respectively, just after
the birth of the fittest clone (t=30 days), when the fittest
population starts to invade (t=36 days), when the tumour returns to
growing almost spherically with low grade of heterogeneity (t=50
days).}

\label{fig0} 
\end{figure}

Considering the dynamics observed in a number of computational
experiments performed at biologically meaningful parameter values, tumour
growth generally appears as exponential, with a consequent
linear increase in the radius of gyration $R_g$ with time{\it , as
  previously observed \cite{Bru19984008}. The main reason, as
  explained in Ref.~\cite{jeon2010off}, is that the growth is driven by those cells that
  reside at the periphery of the mass.} The nearly
spherical shape of the tumour when only a single clone is present
changes significantly in the presence of mutations. The
greatest contribution to asymmetry occurs when a new population with
a faster proliferation rate than neighbouring cells is generated. In that case,
this new population forms an evolutionary niche that can alter the
sphericity of the tumour, until the new clones have proliferated enough
to surround the slower cells and recreate a spherical appearance, as
shown in Fig. \ref{fig0}. Note that, in our model, cells acquire a new phenotype upon mutation in a purely stochastic way and there is equal probability to mutate from the \textit{original} phenotypes to all the others.

The heterogeneity of the mass increases with mutations
until the faster clones are not outnumbering the other
phenotypes. If this occurs, then the Shannon index
$H$ rapidly decreases with time in a way that it is inversely proportional to
the growth of the more proliferating clones, i.e. the faster they grow
the faster $H$ decreases.
The action of the immune system usually tends to favor homogeneity
over heterogeneity, rebalancing the distribution of phenotypes as long
as the immune response is active. As T-cells erode the tumour, 
natural selection leads to an evolutionary bottleneck
characterised by low $H$. It is interesting to note that the roughness
of the tumour tends to be in the interval $1 < M \leq 1.5$, with signs
of superficial infiltration by the
immune system. The limited life span of the CTLs used in our model tends
to promote attacks that occur on the periphery and rarely result in 
deep infiltration, which, as we will show shortly, is instead present when therapies are activated. 



\subsection*{The effect of therapies on tumours}

In the following, all the parameters have been set as
in Ref.~\cite{kim2012modeling} and are reported in Tables
\ref{tab:param_dde} and \ref{tab:param_abm}. The tumour mutation rate has been chosen following
the principle that positive mutations, i.e. mutations that lead to an
evolutionary benefit to the cells over therapies and immune response,
are rare. As it is expected, the dominant phenotype
population usually appears after a clonal expansion of few 
mutations. The range of variation of the proliferating, aggressive tumour has been
set to ensure a biological meaning and a rate of growth
that allows the cancer mass to escape the control of the immune
system in a limited range of time. The main reason for this choice is that we are interested in modeling the impact of different therapies on cancers that will not
be eradicated in the absence of anti-tumour therapies and that, at the same time, do not show growth rates that are unrealistic. Thanks to the probabilistic structure of the system, simulations can generate different outcomes also when parameters are kept fixed for the particular cancer studied. 
Among the different experiments, three paradigmatic dynamics
emerge, which bear particular relevance and help understanding the
typical scenarios that our model predicts. They are the result of stochastic variations on the evolution of initially identical tumours. These cases respectively correspond to a tumour mass
with an initial slow growth and high heterogeneity (case A), and two fast
growing tumours with either initial low (case B) or high (case C)
heterogeneity. 

\subsubsection*{First single-therapy strategy: chemotherapy}\label{res:chemo}
Our first choice is to simulate a cytotoxic chemotherapy that acts with more efficiency against the
cancer cells that have the largest growth rate, starting at
day $60$ after the first tumour cell colonizes the site and for a total
duration of $10$ days. 
The probability of a cell to be killed by chemotherapy, with a total dose of drug  fixed, is set independently from the time duration
of the protocol in an effort to maximize the efficiency of the
therapy.

\begin{figure}[htbp!]
\centering\noindent
\makebox[\textwidth]{
\begin{subfigure}[]{0.68\textwidth}
\includegraphics[width=1\textwidth]{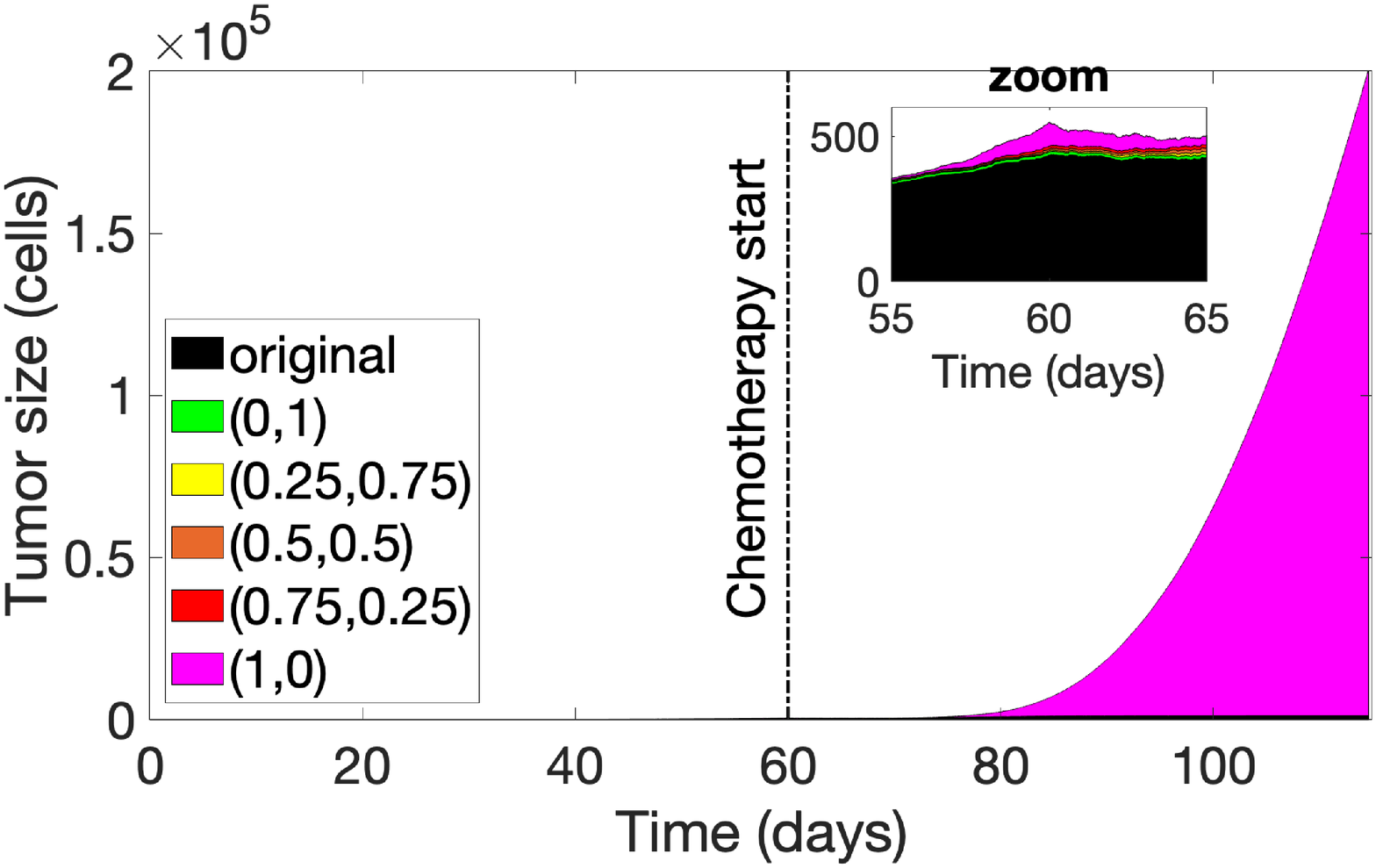}
	\caption{Case A}
\end{subfigure}
\begin{subfigure}[]{0.68\textwidth}
\includegraphics[width=1\textwidth]{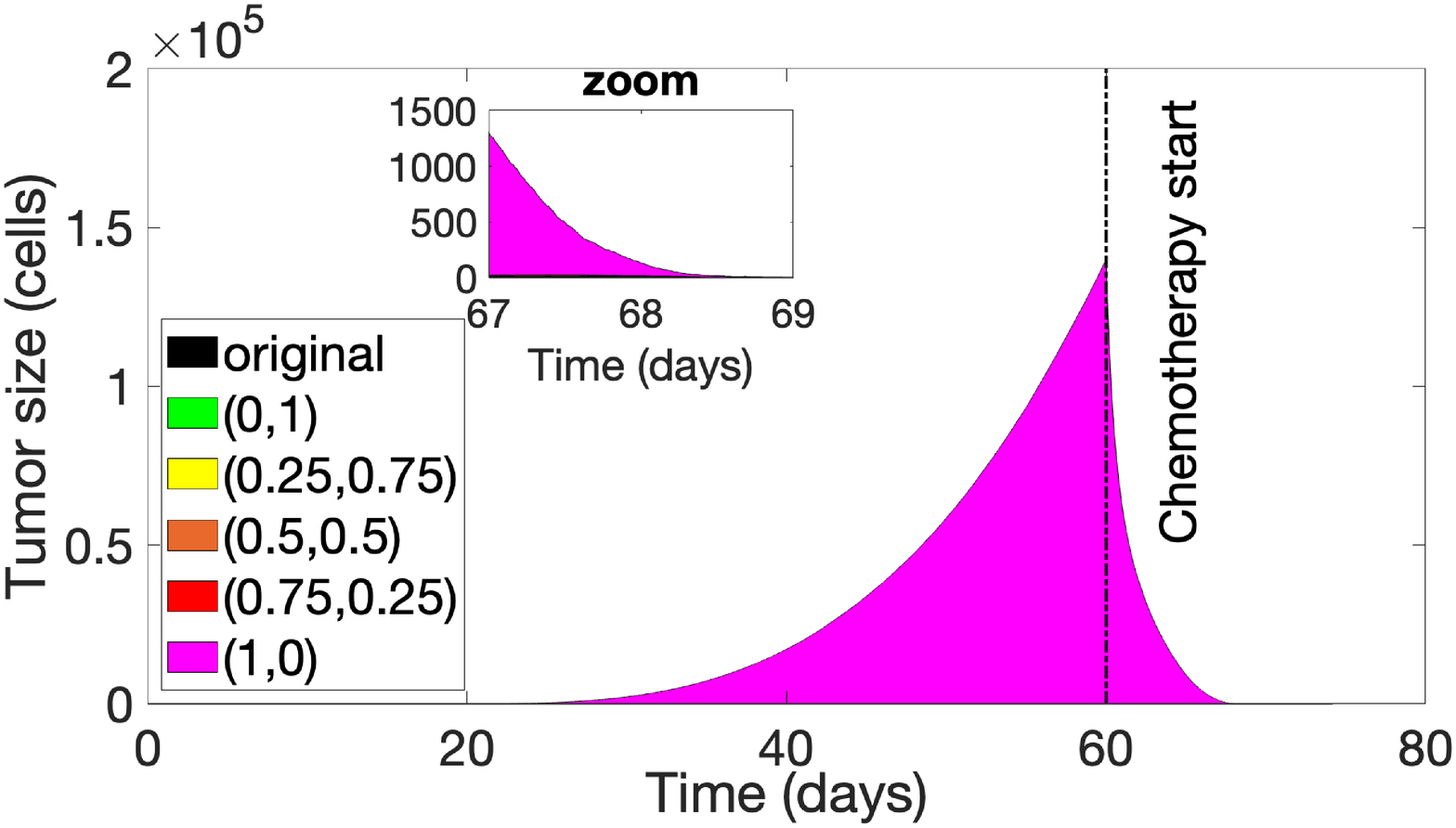}
	\caption{Case B}
\end{subfigure}}\\

\begin{subfigure}[]{0.68\textwidth}
\includegraphics[width=1\textwidth]{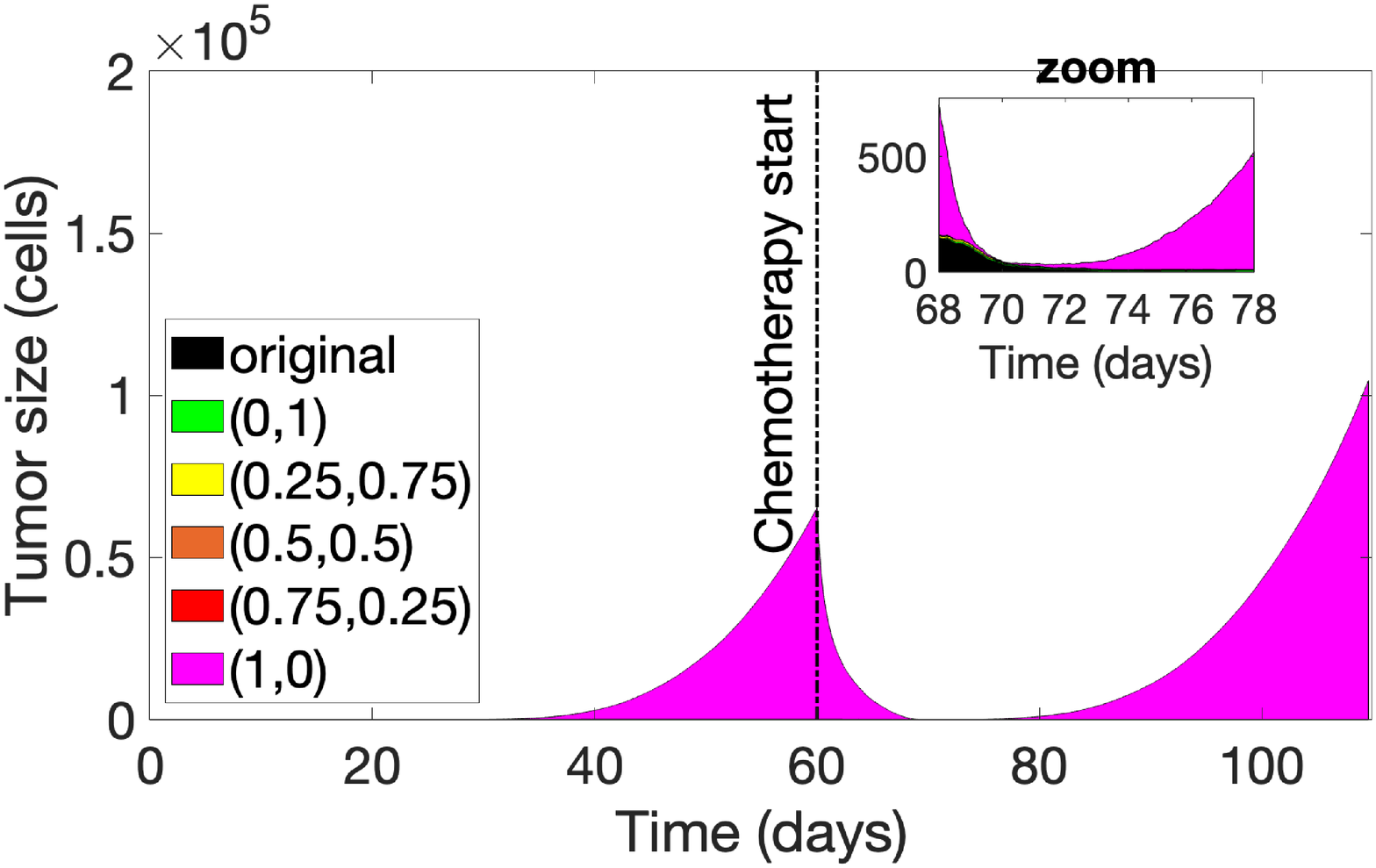}
	\caption{Case C}
\end{subfigure}
\caption[Chemotherapy effect]{{\bf Chemotherapy effect}
Number of cancer cells and phenotypical composition as functions of time. Chemotherapy starts at
day $60$ and lasts for $10$ days. Panel (A) shows an initial phenotype
(in black) with slow reproductive rate (case A), eventually overtaken by a new
phenotype that reproduces very fast (in pink). The inset shows 
cellular distribution during the last part of the treatment. Panel (B)
shows an initially fast reproducing tumour with low heterogeneity (case B),
with an inset of the last two days before complete
eradication. Panel (C) depicts an initially fast reproducing tumour
with high heterogeneity (case C), with the inset focusing on the dynamics
at the end of the treatment and in the first days of the rebound phase.}
\label{au:fig1}
\end{figure}

Fig.~\ref{au:fig1} shows the evolution of the tumour cell population as a function
of time, according to different
phenotypic compositions. In the following, {\it time is evaluated
  starting from the instant at which the original clone starts to
  colonize the new organ. This time does not refer to the primary
  tumour or the history of his evolution.}. In 
panel (A), the effect of chemotherapy on tumour case A, which is representative of those
cancers with lower rates of proliferation but higher
propensities to mutate. 

A too early start of the treatment results 
in a completely ineffective strategy, with a negative outcome. This is
because chemotherapy reduces the more proliferating cells (in pink) at day $60$
when those cells are still scarce and the tumour is too small to
benefit from the action of the cytotoxic drug. Once the treatment is
over, the remaining cells restart to mutate and proliferate, with an
exponential growth that the immune system alone cannot contain. As
shown in the inset of panel (A), the number of cells belonging to the original phenotype (in black) remain almost
constant throughout the procedure and do not change significantly for the duration of the
experiment. 

Panel (B) of Fig.~\ref{au:fig1} instead shows a complete eradication of
case B, where tumour cells have initially a low heterogeneity but are
reproducing fast. The effect of the therapy 
is in this case to eliminate every cell belonging to the dominating,
fast-reproducing phenotype before it is over, i.e. approximately
around day $7$ of its $10$ days duration. Also, all cells
of the original phenotype are eradicated by the end of the treatment,
with the tumour completely cleared out by the effect of the cytotoxic
drug and the immune response.

An initially fast reproducing tumour with high heterogeneity can instead
lead to uncontrollable rebounds, with an overall negative outcome for
the patient. In panel (C), the action of the chemotherapy is not
sufficient to eliminate every single cell belonging to the mutated
phenotype. According to our choice of parameters, it is enough that one original
cancer cell or one of the more proliferating clone survives after the
chemotherapy that a fast, uncontrollable rebound can be expected.

Interestingly, these last two cases (i.e. B and C) do not show different evolutions
of the radius of gyration $R_g$ {\it (not shown in the Figure)}
during the action of chemotherapy, since the treatment acts homogeneously
on the cancer mass as a whole. This is associated directly to the 
limited dimension of the tumour, leading to the drug acting
on the aggregate with a high strength from all spatial directions. 
Roughness $M$ instead show significant
changes from case A and cases B and C. Tumour case A 
remains spherical and compact during the
experiment, essentially because the treatment has a very limited
effect on the mass due to its premature start. Tumour case B and C, instead, reach high
level of roughness during the treatment, showing, for the cases
reported in the Figure~\ref{au:fig1}, a maximal $M$ of $3.77$ and $3.45$,
respectively at day $64$ and $65$ for cases B and C. This clearly
indicates that the tumour
loses density and becomes morphologically inhomogeneous at around half
of the treatment duration and is infiltrated to a relevant degree by
the T-cells taking part in the immune response. The overall
consideration from these results is that the correct timing of
treatment, {\it here intended as the ideal treatment
starting time and therapy duration to achieve optimal patient’s benefit}, 
is a major variable for the outcome of the treatment and
it is also strongly affected by phenotypical compositions. 

\subsubsection*{Second single-therapy strategy: immune boosting}\label{res:immuno}

\begin{figure}[htbp!]
\centering\noindent
\makebox[\textwidth]{
\begin{subfigure}[]{0.55\textwidth}
\includegraphics[width=1\textwidth]{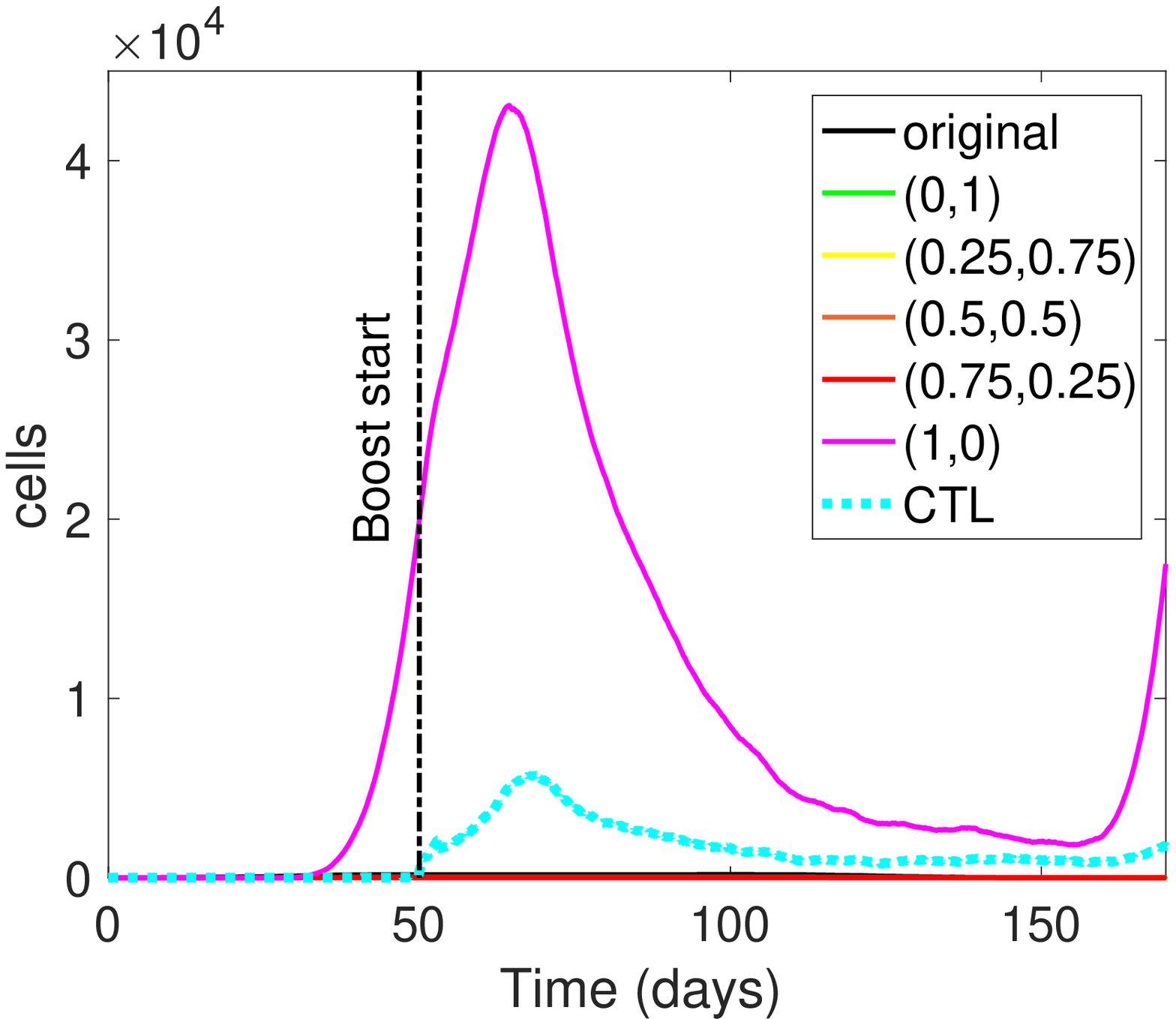}
	\caption{}
\end{subfigure}
\begin{subfigure}[$t=170 \mathrm{days}$]{0.59\textwidth}
\includegraphics[width=1\textwidth]{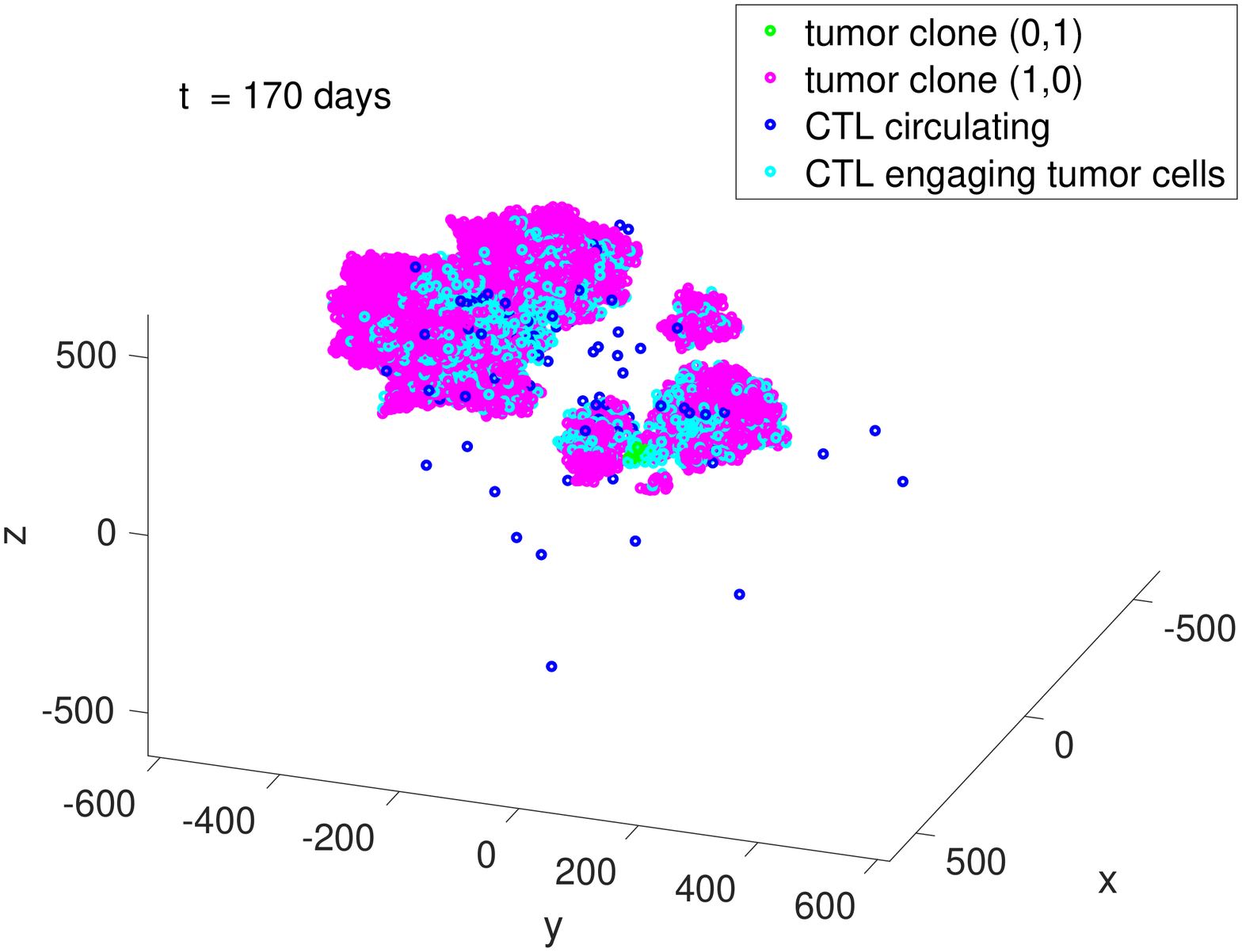}
	\caption{$t=170 \mathrm{days}$}
\end{subfigure}}
\caption[Immunotherapy effect.]{{\bf Immunotherapy effect.}
Results for the prototypical case C, a fast growing cancer with high
phenotypical heterogeneity. Boost immunotherapy starts at day $50$ and
lasts a total of $3$ days. Panel (A): time series of each cell
population, panel (B): snapshot of the tumour mass at time
$t = 170$ days.} 
\label{au:fig2}
\end{figure}

The overall effect of immune boosting is to increase the number of CTL cells circulating
around the tumour site, which we simulate as an injection of cells
starting at day $50$ 
and occurring for a duration of $3$ days. Erosion of cancer cells by the immune system proceeds
from the periphery towards the center of the tumour mass, and is
usually characterised by a linear decrease of $R_g$ during the first
few days. Another typical characteristics of the dynamics that
follows boosting is a clonal expansion of the CTL
population shortly after treatment. For the prototypical three
cases A, B and C introduced above, all of our computational experiments
indicate that boosting alone is not able to
eradicate cancer: after an initial decrease in the tumour mass, two types of evolutions have been
observed, both resulting in negative outcomes. Of particular relevance
is case C, which, although not treatable by chemotherapy alone, shows a somewhat unexpected and
complicated morphology when
subject to immune boosting. In fact,
after an initial clonal selection of the less immunogenic phenotype, 
case C displays a clear deviation from sphericity in the mass, with a
nonlocal spread of the tumour population in islands of different sizes,
as reported in Fig.~\ref{au:fig2}. After a decreasing phase
due to an immune response that does not result in a complete
eradication, the tumour population is eventually subject to
a faster, uncontrollable increase driven by disconnected, smaller masses.
Overall, a selection of 5 parameter sets and 10 trials for different seeds give qualitatively similar results. 

Panel (A) in Fig.~\ref{au:fig2}, shows the number of cells for each
phenotype as a function of time: an initially exponential growth
is firstly slowed down then halted by immunotherapy, with a maximum cell population
occurring about $20$ days after the beginning of the
treatment. Around day $100$, the evolution of the cell population
changes. Frequent, local maxima in both tumour and T-cell curves represent
the failed attempts made by the immune system to completely erode the tumour due to the increasing sparseness of the
cancer. This behaviour seems to occur for a protracted period of time
of about $60$ days.

As the cells in the island's sizes begin to proliferate faster than
the rate of killing of the T-cells, a rebound phase with a higher
speed of growth than the \textit{original} unbroken mass appears at day $180$. Panel (B) provides an image
of the scattered status of the tumour immediately before its exponential
rebound. Let us remark that the model does not allow for migration of cancer cells and this picture is the result of the 
infiltration of T-cells coming from boosting and immune
response. 

As expected, morphology immediately before the rebound phase
is characterised by a high value of roughness, with
$\max(M_{\mathrm{Case C}})=2.16$ at day $162$. Also, there is almost a
twofold increase in $R_g$ compared to the value for identical number
of cells in the first growth phase. For example, for $10^4$ cells, we have $R_{g_{\mathrm{Case
      C}}} = 125$ at day $47$ 
and $R_{g_{\mathrm{Case C}}} = 240$ at day 160. 

\subsubsection*{First example of synergistic therapy: chemotherapy and
boosting}\label{res:combi}
Cancer heterogeneity has been invoked to explain one of the major
aspects of cancer development, namely acquired drug resistance, by
which phases of remission are often followed by a rapid growth of
tumour cells~\cite{gerlinger2010darwinian}.
One of the ways to overcome resistance is, for instance, to find more ``evolutionarily
enlightened'' strategies that places malignant cells in an
``evolutionary double bind''~\cite{gatenby2009lessons}. In cancer, a double blind could be
obtained using the immune system as natural biological
predator~\cite{basanta2012exploiting}. Clinical evidence shows that
immunotherapy or oncolytic viruses alone are not effective, despite
the possible theoretical advantages. Therefore, cancer treatment is
adopting a multistep approach that combines biological and
chemical$/$radioactive therapies using cytotoxic effects on one side
and subsequent adaptation on the other side to limit tumour adaptive
resistance~\cite{ramakrishnan2010chemotherapy,
  antonia2006combination}.

\begin{figure}[htbp!]
\centering\noindent
\makebox[\textwidth]{
\begin{subfigure}[Case A]{0.65\textwidth}
\includegraphics[width=1\textwidth]{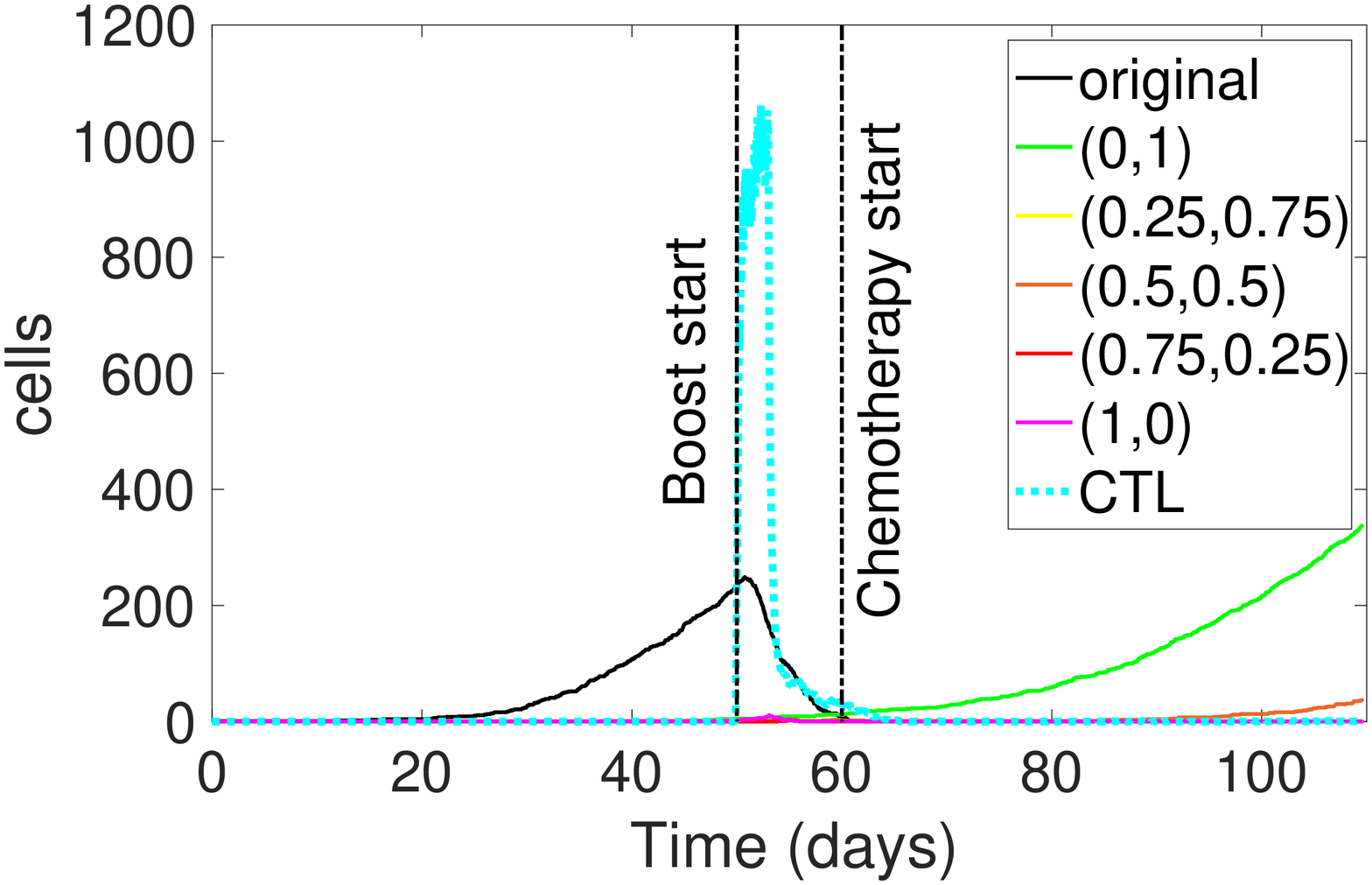}
	\caption{Case A}
\end{subfigure}

\begin{subfigure}[Case B]{0.65\textwidth}
\includegraphics[width=1\textwidth]{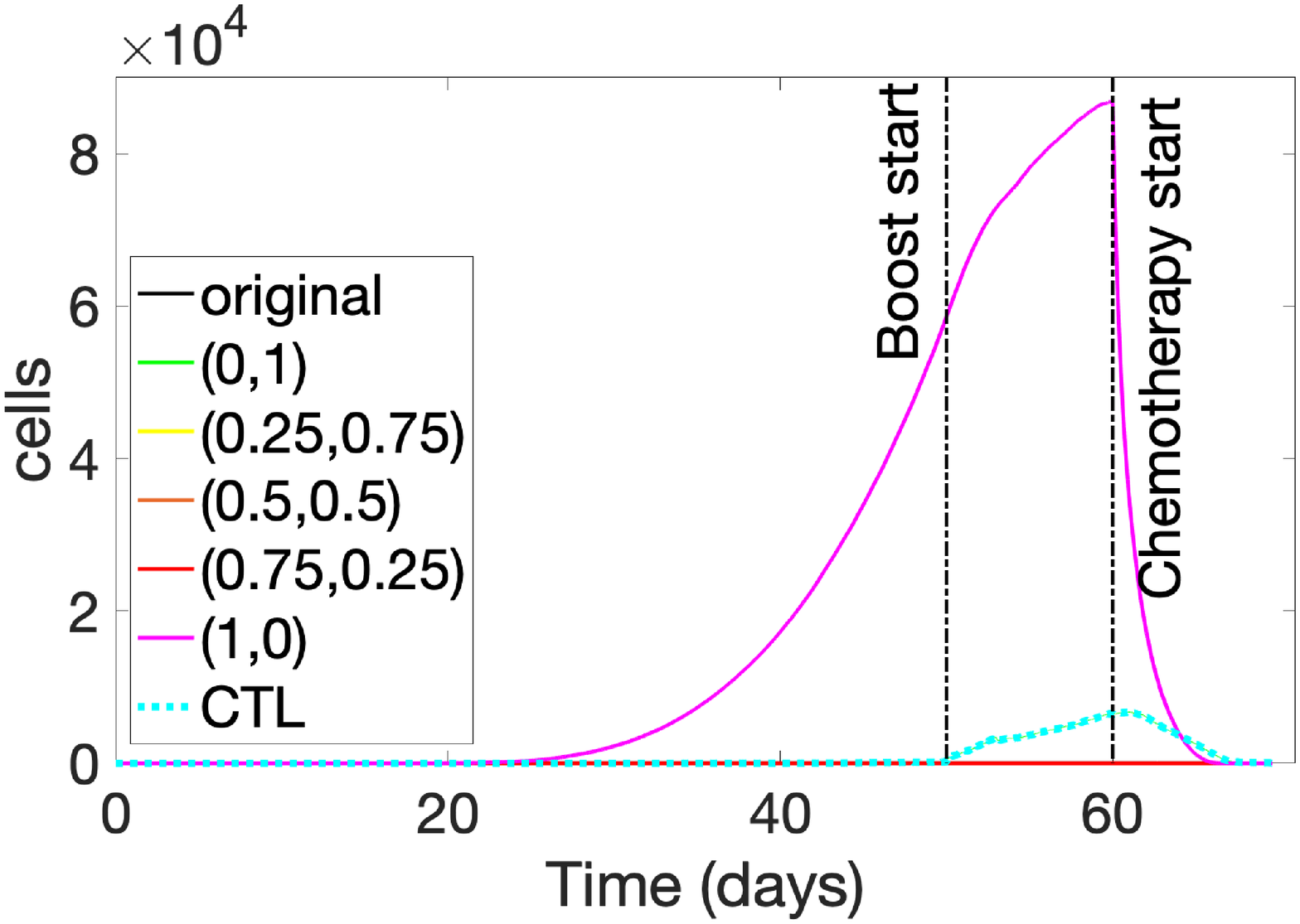}
	\caption{Case B}
\end{subfigure}}\\

\begin{subfigure}[Case C]{0.65\textwidth}
\includegraphics[width=1\textwidth]{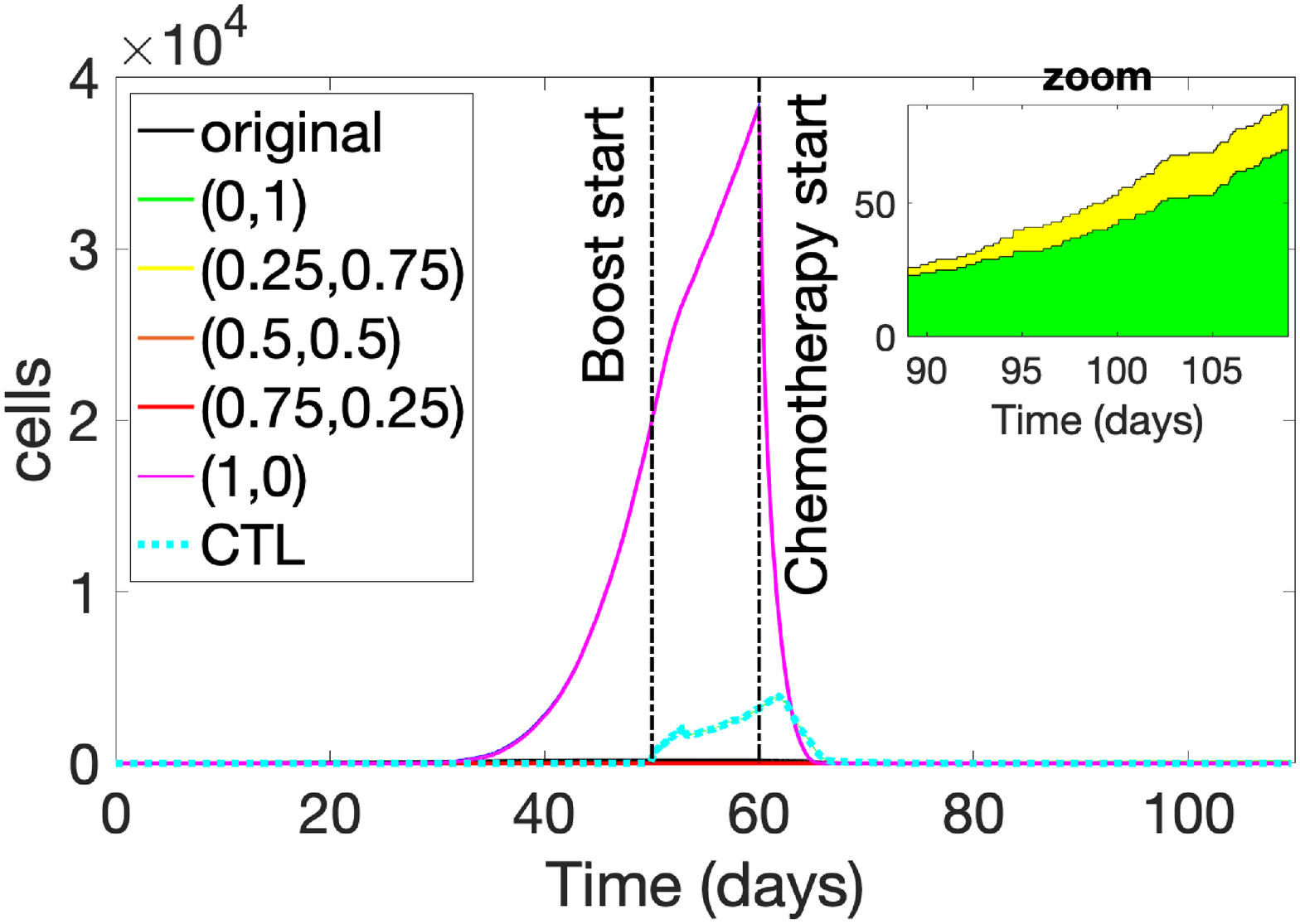}
	\caption{Case C}
\end{subfigure}
\caption[First kind of synergistic therapy: chemotherapy and
    immune boosting.]{{\bf First kind of synergistic therapy: chemotherapy and
    immune boosting.}
Plots represent the time evolution of cellular populations, with
insets showing phenotypical composition of the tumours with time.
Immune boosting 
starts at day $50$ and lasts $3$ days, whereas chemotherapy starts at
day $60$ and has a duration of $10$ days. Panel (A), (B) and (C)
respectively show cases A (slowly growing tumour), B (quickly growing tumour, low
heterogeneity) and C (quickly growing tumour, high
heterogeneity). {\it Note that different scales have been used to allow for
greater details of the dynamics. The
inset in Panel (C) represents a close-up of the time range $90-110$
days. The unit on the $y-$ axes of every plot in the Figure is cells' numbers.}} 
\label{au:fig3}
\end{figure}

Guided by the poor outcomes displayed by immune boosting alone in the
prototypical cases, we now consider the
combination of chemotherapy and immune boosting, with the aim of
discussing the major factors that maximize positive outcomes. The
prototypical cases have been subjected to a protocol of an immune boosting
injection at day $50$ lasting three days, followed by a chemotherapy
session at day $60$. Results are displayed in Fig~\ref{au:fig3}, with the
insets displaying phenotypical composition over time. 
For cases A and B,
the complete temporal range is shown, whereas for C the last $20$
days are reported. Timing for these therapies has been chosen
arbitrarily. {\it For cases B and C the second lesion grows up
to numbers of tumour cells that are close to
the detectability threshold.} 

The effects of this synergistic therapy in cases A and C are
similar: chemotherapy preferentially kills those cells that are
fast to reproduce, leaving the slowest reproducing phenotype unaffected. As a
result, rebounds occur once therapies end, with case A showing a
negative outcome within the simulated time window and case C
displaying a still moderate but uncontrollable growth at the end of
the simulation. In other words, the effect of chemotherapy is
to create an evolutionary bottleneck
that selects the poorly immunogenic clones. {\it In particular, case C (see
the inset of Panel (C)) shows a surviving tumour composed by only two
clones: the clone $(0,1)$ and the clone $(0.25,0.75)$: these are the two
families that are the slowest in proliferating and have the lowest
immunogenicity.}
These clones have a strong
immunoediting ability and remain unnoticed by the immune system,
resulting in unending growth: {\it because of their ability to be
  elusive to both the immune response and chemtherapy, outcome C
  appears the worst of all.} Note that, for this reason, this phase
is of a different nature than those previously reported for individual therapies
(i.e. Figs.~\ref{au:fig1} and ~\ref{au:fig2}). Also, because of the low
heterogeneity of case B, this combined therapy is instead successful
in fully eradicating the tumour, which is eliminated during the
administration of the cytotoxic drug. On the other hand, case A, differently from the
results obrained for chemotherapy alone, displays a selection of the
poorly differentiated immunogenic clones. 

From the morphological perspective, masses emerging from this
synergistic intervention appear to be low in roughness when the
reproductive rate is slow, with
$\max(M_{\mathrm{Case A}}) = 1.62$ at $t = 81$ days. If the rate is
instead fast, the level of heterogeneity usually determines the level
of roughness, with
low heterogeneity contributing high roughness during the chemotherapy
phase, i.e. $\max(M_{\mathrm{Case B}}) = 3.05$ at
$t=64$ days and $\max(M_{\mathrm{Case C}})=2.64$ at $t=63$ days.
Contrary to case B, case C shows a Shannon index of $H>0.5$ for most
of the simulation, which results in an unsuccessful eradication. 

Other time protocols and order of administration between boosting and
chemotherapy are possible, and have been tested to some degree
(results not
shown here). Although a study of optimisation of protocols is not
within the scope of the present work, the overall insight from the
simulations is that heterogeneity always plays an important role in the
outcomes. For this combination therapy, high values of $H$ are 
consistently associated with negative prognosis \cite{greaves2012clonal}.

\subsection*{Second example of synergistic therapy: radiotherapy,
boosting and ``abscopal'' effect}\label{res:abs}

A second example of synergistic therapies that is currently used in
clinical practice is the combination of an initial cycle of
radiotherapy with an immune boosting protocol. Besides a better
understanding of the parameters that can trigger a positive outcome,
our interest in testing such a combination resides in the occurrence
of a somewhat rare and poorly understood
event, which is named ``abscopal'' effect. There are a number of clinical cases discussed in the medical
literature where a reduction of a
secondary tumour or an existing metastasis outside the primary, radiated
lesion has been observed \cite{reits2006radiation, lugade2005local,
  finkelstein2012clinical, finkelstein2012combination,
  finkelstein2011confluence, vatner2014combinations}. Differently from
chemotherapy, radiotherapy has a localised action on the region
irradiated and this makes the phenomenon, to some extent,
counterintuitive. 
Sometimes, the effect appears to affect a secondary lesion very distant from the region treated.

The
complications inherent to the stimulation of such an event are due
to the immune action apparently being as crucial as radiotherapy in triggering such a response. 

The model allows us to test some hypotheses on the nature and causes
of the effect of protocols introduced by Demaria et al. in
Ref.~\cite{demaria2004ionizing}, who have reported some interesting
and positive outcomes. In particular, they have treated mice with a
syngeneic mammary carcinoma in both flanks with immunotherapy and only
one of the two tumours with radiotherapy. 
They use the non-irradiated lesion to monitor the insurgence of the abscopal effect. 
Biologically, a tumour-specific T-cells activation occurs after 
inflammatory signals are introduced in the system as a consequence of
the therapy. Dying cancer cells
release tumour-specific antigens and immune-stimulatory signals that seem to
induce an increased 
recognition of cancer cells with the same phenotypical
characteristics in other areas of the body. Further, radiation modulates different
compartments of the tumour microenvironment, resulting in
exclusion-inhibition of effector T-cell and induction of de novo
anti-tumour immune responses~\cite{demaria2015role}. The protocol that
we simulate is a radiotherapy (RT) on the primary tumour (not
simulated or showed here) at day 1,
followed by an immune boosting that lasts 10 days. 
As anticipated, RT is considered a restoring factor in the ability
of CTL cells to recognise and kill all cancer phenotypes, with no
exceptions.

\begin{figure}[htbp]
\centering\noindent
\makebox[\textwidth]{
\begin{subfigure}[]{0.50\textwidth}
\includegraphics[width=1\textwidth]{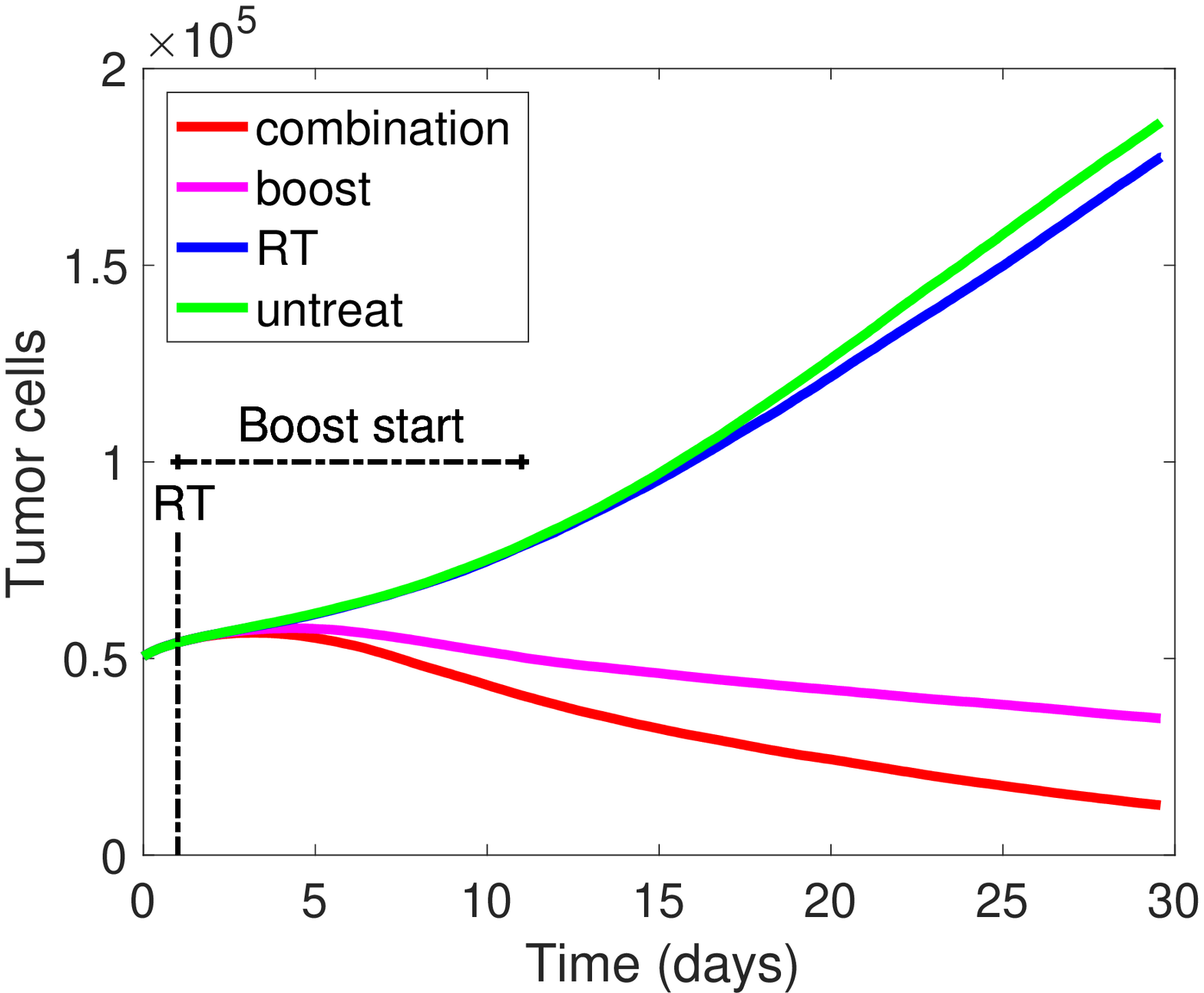}
	\caption{Therapy's action on tumour}
\end{subfigure}
\begin{subfigure}[]{0.50\textwidth}
\includegraphics[width=1\textwidth]{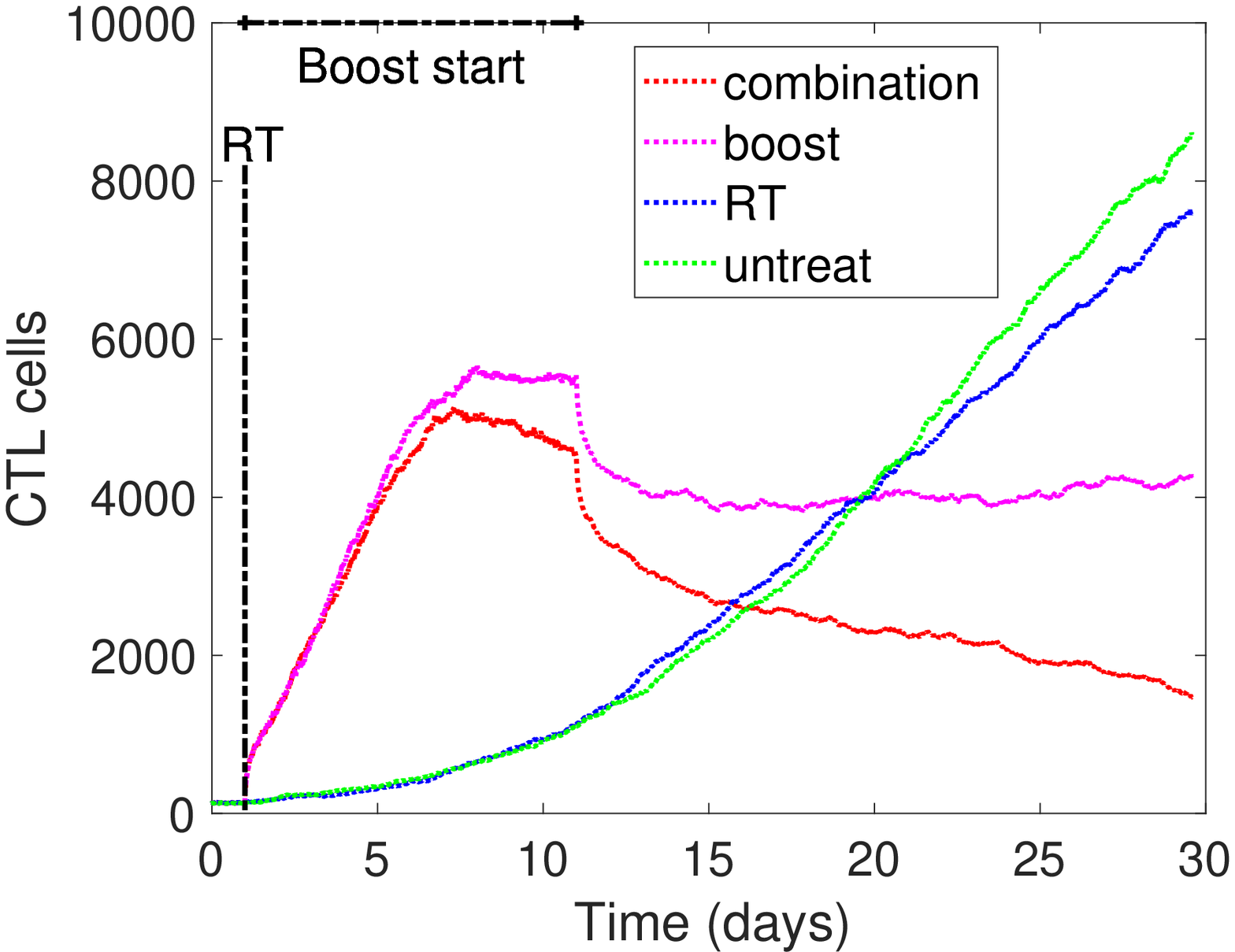}
	\caption{Therapy's action on CTL}
\end{subfigure}}\\

\centering\noindent
\makebox[\textwidth]{\begin{subfigure}[]{0.50\textwidth}
\includegraphics[width=1\textwidth]{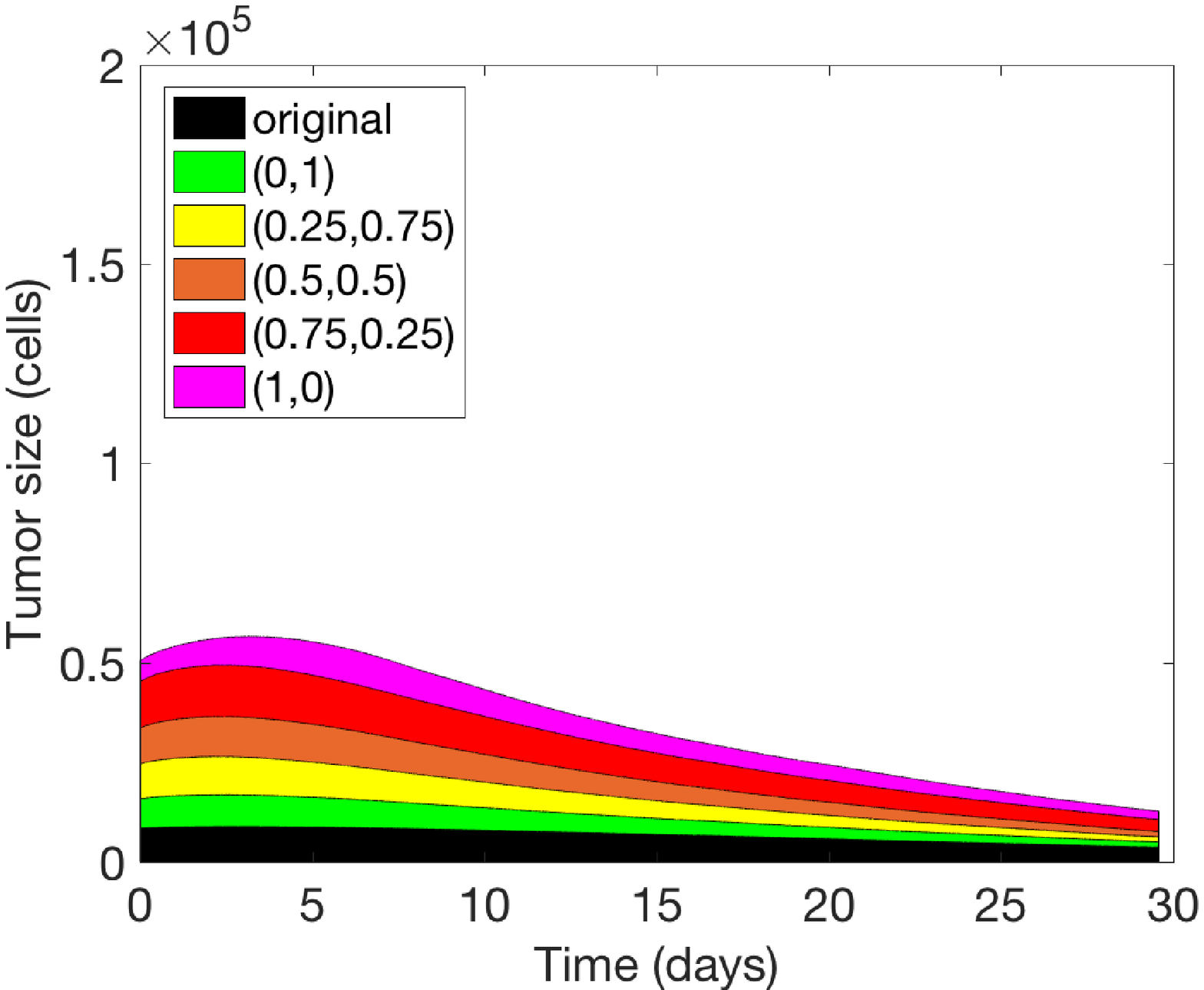}
	\caption{Combination}
\end{subfigure}
\begin{subfigure}[]{0.50\textwidth}
\includegraphics[width=1\textwidth]{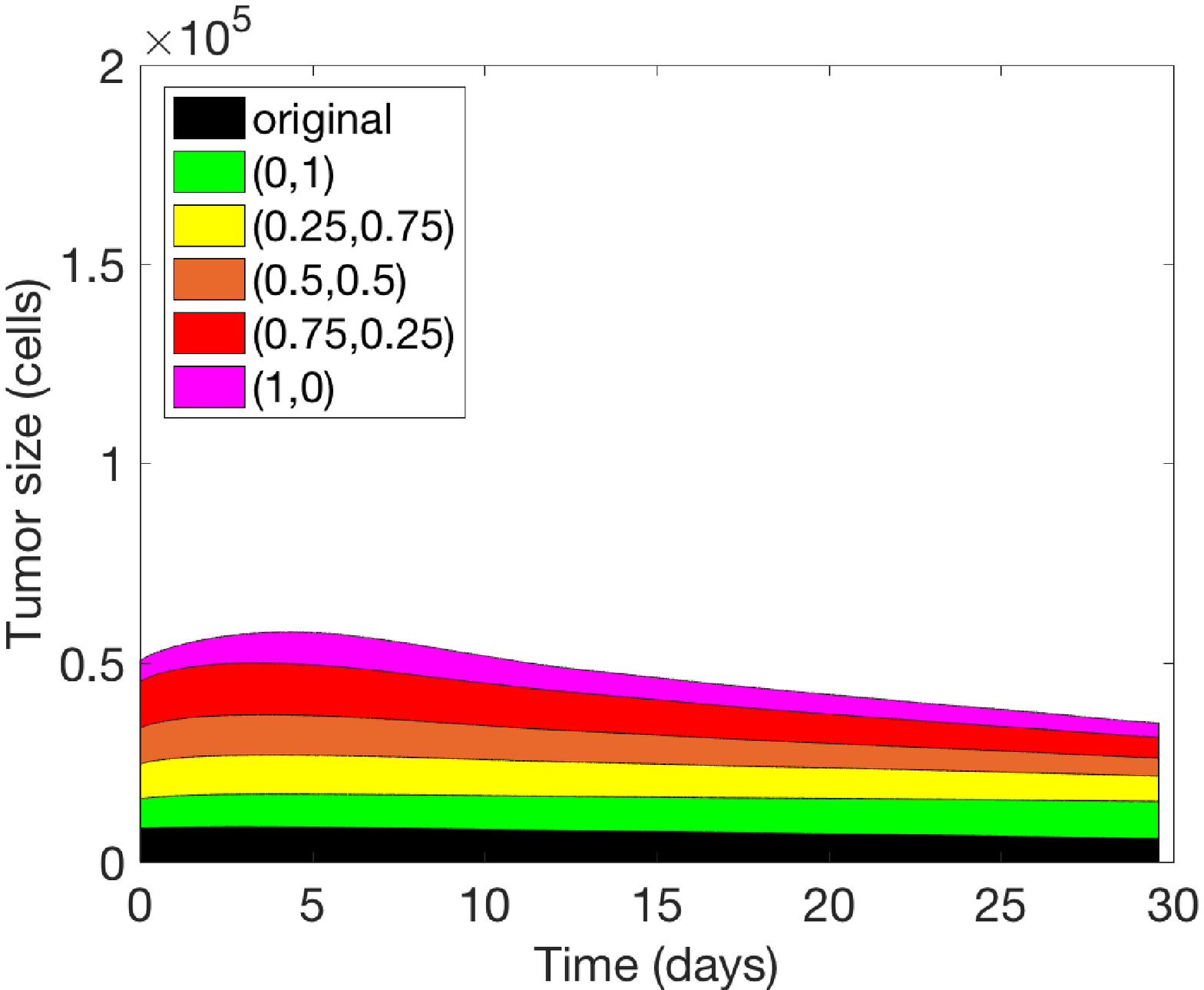}
	\caption{Immunotherapy}
\end{subfigure}}\\

\centering\noindent
\makebox[\textwidth]{\begin{subfigure}[]{0.50\textwidth}
\includegraphics[width=1\textwidth]{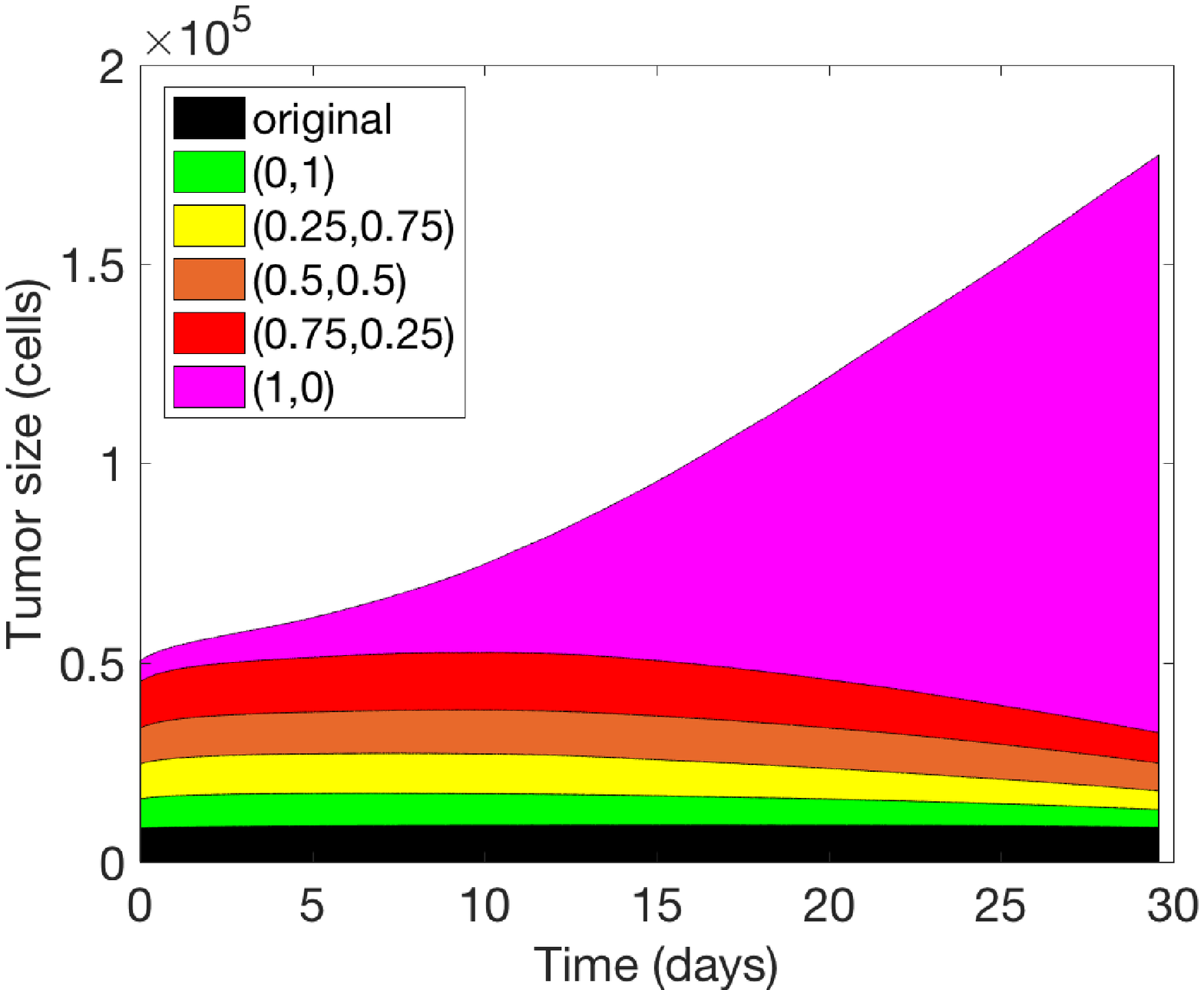}
	\caption{RT}
\end{subfigure}
\begin{subfigure}[]{0.50\textwidth}
\includegraphics[width=1\textwidth]{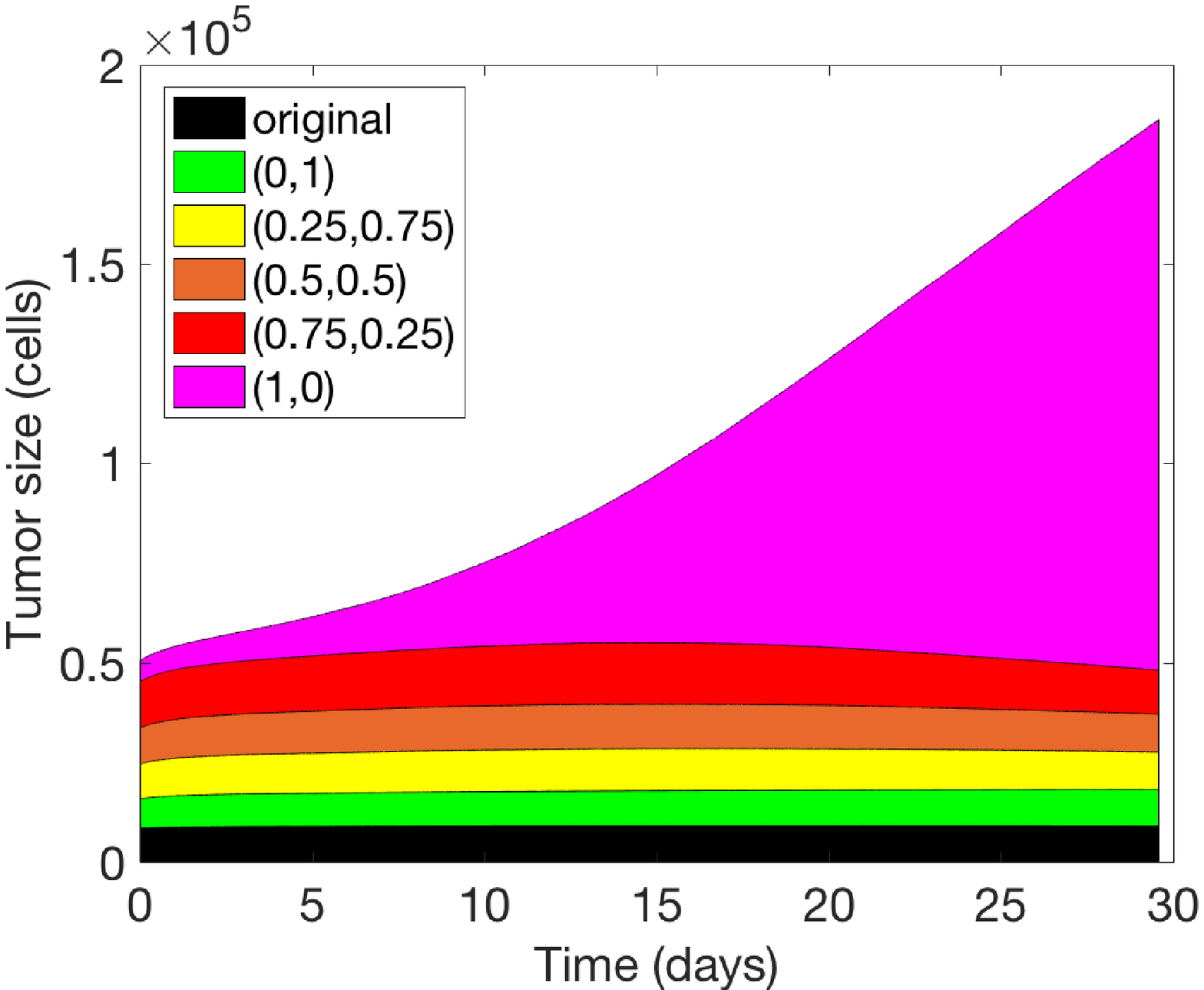}
	\caption{Control case}
\end{subfigure}}
\caption[Effects of anti-cancer treatments on a secondary tumour.]{{\bf Effects of anti-cancer treatments on a secondary tumour.}
Panel (A and B): Cellular populations of a
secondary tumour as functions of time, with cancer populations depicted
by a continuous line (panel (A)) and CTLs population by a dashed line
(panel (B)). Each color corresponds to a different treatment. {\it These curves represent
 a typical simulation run. Variability due to stochasticity is
 minimal.} Other panels show
phenotypical composition of the secondary tumour for the {\it curves
  shown in Panels (A)-(B)) and for the} following
strategies applied the primary tumour: combination therapy (panel (C)),
immune boosting (panel (D)), radiotherapy (panel (E)), no treatment
(panel (F)). Each color indicates a different cancer clone
family. {\it Legends show labels for such families, whose properties
  are described in
  detail in Table~\ref{tab:clones}.}}
\label{au:fig4} 
\end{figure}

The secondary lesion is composed of $5$ x $10^4$
heterogeneous cells (the same number of cells injected in mice in the
experiment by Demaria et al.), 
generated randomly with each clonal family representing at least $10\%$ of the total population.
We compare the
action of two single therapies (RT and immune boosting),
with a combination of the two and a control case where the second
lesion remains untreated. Results are presented in
Fig~\ref{au:fig4}: 
{\it each panel represents the typical outcome from
  a single simulation. For each case, i.e. control (no treatment),
  combination, immunotherapy, and radiotherapy, we have performed three
  different simulation runs, with different initial conditions. The
  outcome of each simulation for any configuration consistently gives
  comparable results. Variation due to stochasticity are minimal and
  do not affect the outcomes.}

Firstly, no treatment or RT alone result in similar negative outcomes for the
secondary mass, not directly treated by RT, both
from the perspective of surviving cancer cells (panel (A)) and the
response of the CTLs of the immune system (panel (B)). An initial RT
with no follow-up has the only effect of delaying an exponential
rebound, not dissimilar to the behaviour of an untreated
mass. Boosting alone does not impact the mass as much as when
we combine boosting and RT, with the former giving rise to a tumour
that after $30$ days has less than half the mass than in the case of
chemotherapy alone. Stimulated by the release of the antigens of
the dead tumour cells of the primary lesion, both therapies show a
maximum in the number of active CTLs, which occurs around day $7$ and is then
followed by a characteristic drop around day $10$. Qualitatively the results of the model reproduce the experimental data in Ref.~\cite{demaria2015role}, with indications of a start of a remission for the secondary tumour.

Different strategies on the secondary mass also lead to different
clonal compositions, which have an effect on the final outcome..One of
the keys for the success of the strategy is to have no phenotype
prevailing over the others, as shown in panel (C) for the combination
of RT and immune boosting and, to a lesser
extent, in panel (D) for immune boosting alone. {\it Note that the width of
the coloured regions in panels (C)-(F) indicates the number of cells
that belong to a specific clone population: the larger the width, the
larger the population.} 
For example, combination therapy
provides a very high Shannon index, $H>0.99$, throughout the whole
duration of the experiment. RT and no treatment show instead
reduced indices, with values at day $30$ 
of $0.43$ and $0.55$ respectively. Interestingly, the RT case appears to be less
heterogeneous than the control case. Overall, it is important to stress
that, for the case of the secondary lesion, high heterogeneity is not
per se associated to a worse prognosis.  
The reason is that a successful
action on the secondary mass 
reduces the fitness advantage of the phenotypes and makes the immune
system able to recognise each clones equally.  {\it Note that Panel (E)
refers to the radiotherapy case made on primary tumour, showing the
evolution of the secondary mass and the fact that the immunogenic
effects induced by the treatment are not sufficient to contrast cancer
growth.  Panel (F) represents the control case, where no treatments
are administered and the tumour is growing unchallenged as an
aggressive breast cancer.}


We find that the relation between the phenotypical composition of the primary and
secondary lesions plays a very important role in the dynamics of the
so-called abscopal effect. The previous examples refer to
a secondary lesion that is antigenically related to the first tumour,
but this is not always the case in practice. Results vary considerably
if the phenotypical compositions differ and this is important to stress.

\begin{figure}[htbp!]
\centering\noindent
\makebox[\textwidth]{
\begin{subfigure}[]{0.55\textwidth}
\includegraphics[width=1\textwidth]{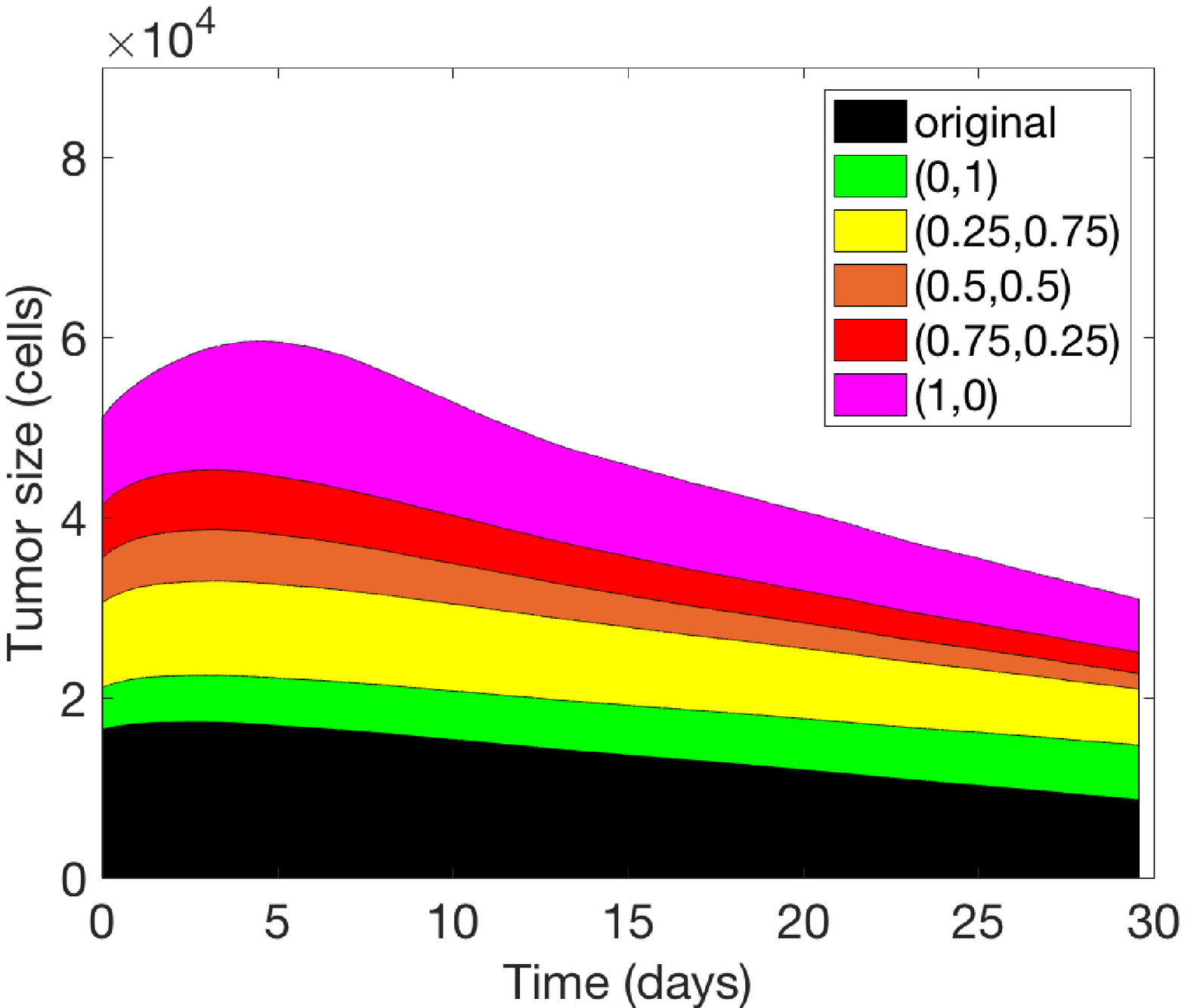}
	\caption{Antigenical unrelated tumour}
\end{subfigure}
\begin{subfigure}[]{0.55\textwidth}
\includegraphics[width=1\textwidth]{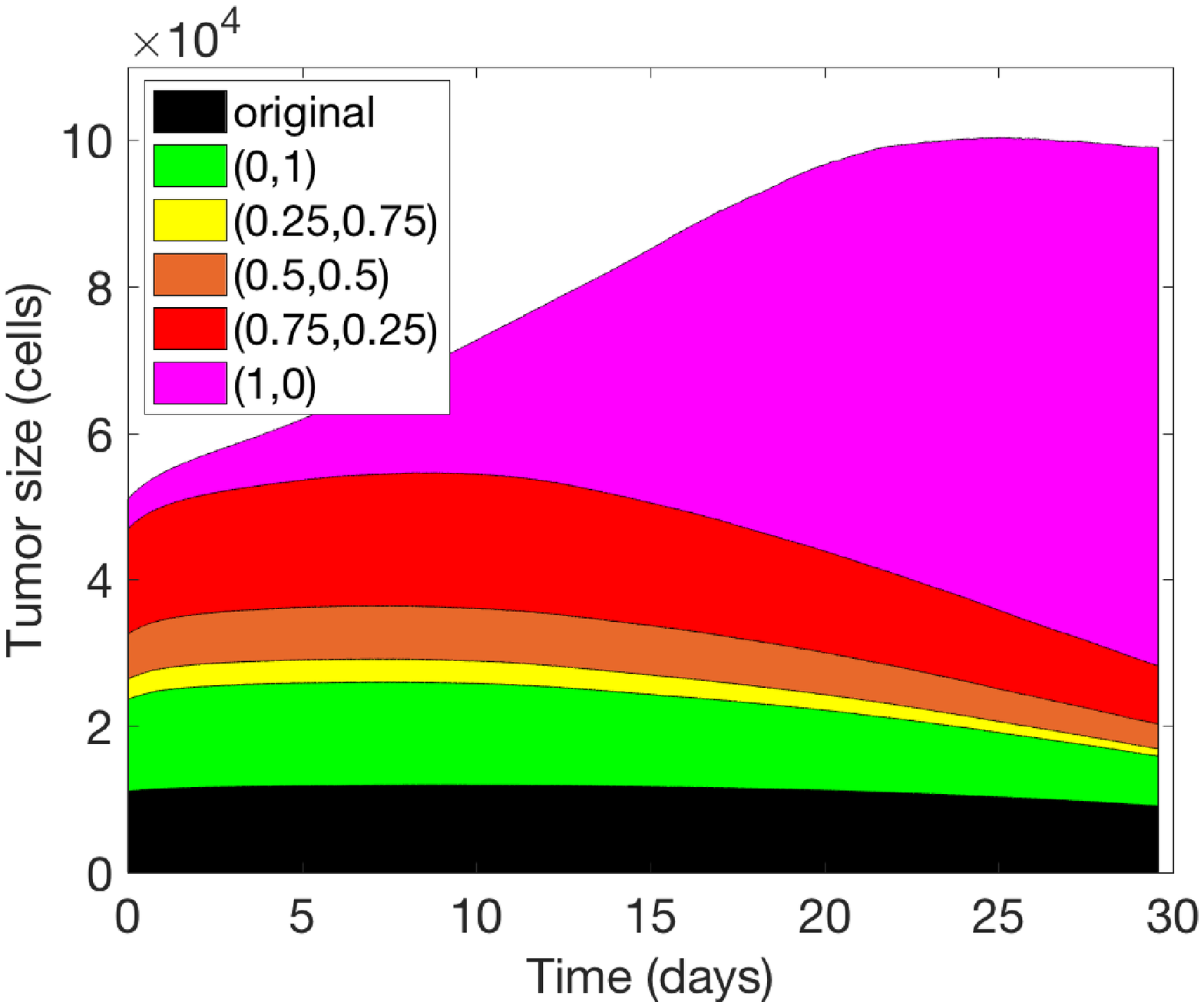}
	\caption{Immune suppressed host}
\end{subfigure}}
\caption[Immune mediation on a secondary lesion.]{{\bf Immune mediation on a secondary lesion.}
Tumour size and composition as functions of time for three different
secondary tumours, subjected to different anti-tumour strategies on the primary
tumour. Panel (A) shows results for a secondary tumour antigenically unrelated to the
primary and panel (B) shows the outcome for a secondary tumour that is antigenically related to the
primary tumour but implanted in an immune-suppressed host.
{\it
  Legends show labels for the tumour subpopulations emerging from the
  “original” population, via secondary mutations caused by the
  therapies. Properties of each clone family are described in detail
  in Table~\ref{tab:clones}.}} 
\label{au:fig5}
\end{figure}

For example, if the secondary tumour is characterised
by clones that are not antigenically related to the first lesion, the
final outcome of combination therapy cannot be as positive as in the
previously discussed cases. Given that radiotherapy affects
phenotypes that are distributed in different ways in the first and
second tumours, the immune system is not capable of recognising specific tumour
cells in the same successful way as in the previous examples. As a
result, the therapy shows a worse outcome, as can be seen in
Fig.~\ref{au:fig5}(a). Furthermore, if the second tumour mass is instead
implanted in an immune suppressed host (mathematically obtained by setting $\alpha=10^{-15}$), 
where the level of CTLs circulating around the tumour is lower than ordinary
levels, the outcome is negative. As reported in
Fig.~\ref{au:fig5}(b), a reduced fitness of the immune system causes
one phenotype to prevail over the others and proliferate quickly out
of control. It is reasonable to suppose that, if more cycles of therapies are repeated, the effectiveness of treatments is likely to be also reduced.
\\

\subsubsection*{Mutation rates and eradication of secondary tumours} \label{res:abs_mut}
A complete eradication of a secondary tumour as
an indirect result of an anti-cancer therapy on a primary lesion is a
relatively rare occurrence. Also, it appears to be associated mainly with certain types of cancer, 
namely melanoma or breast cancer. This
might be linked to the fact that generic metastatic events are
characterised by a high genetic instability, often
making secondary lesions phenotypically unrelated to the first
tumour. From this point of view, a possible speculation could be that
the so-called abscopal effect is not a rare event per se,
but it is an effect limited to secondary tumours that have a
phenotypical clone composition that is not too dissimilar from the
originally treated lesion, and thus the effect only seldom changes the
prognosis for secondary lesions. Indeed, extensive genetic and
phenotypic variation are known to strongly influence therapeutic
outcomes \cite{burrell2013causes}. 

\begin{figure}[htbp!]
\centering\noindent
\makebox[\textwidth]{
\begin{subfigure}[]{0.65\textwidth}
\includegraphics[width=1\textwidth]{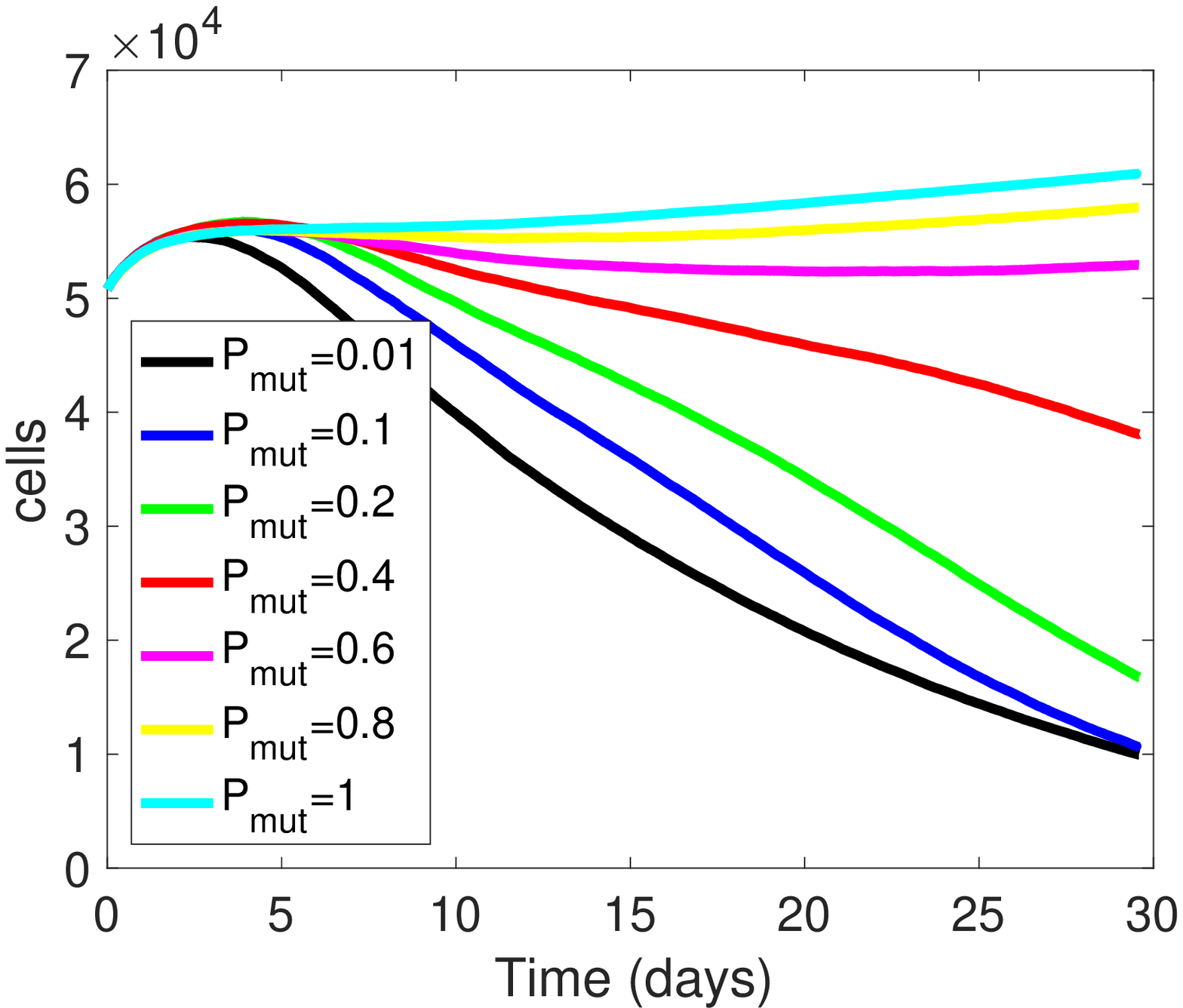}
	\caption{}
\end{subfigure}}\\

\centering\noindent
\makebox[\textwidth]{
\begin{subfigure}[]{0.55\textwidth}
\includegraphics[width=1\textwidth]{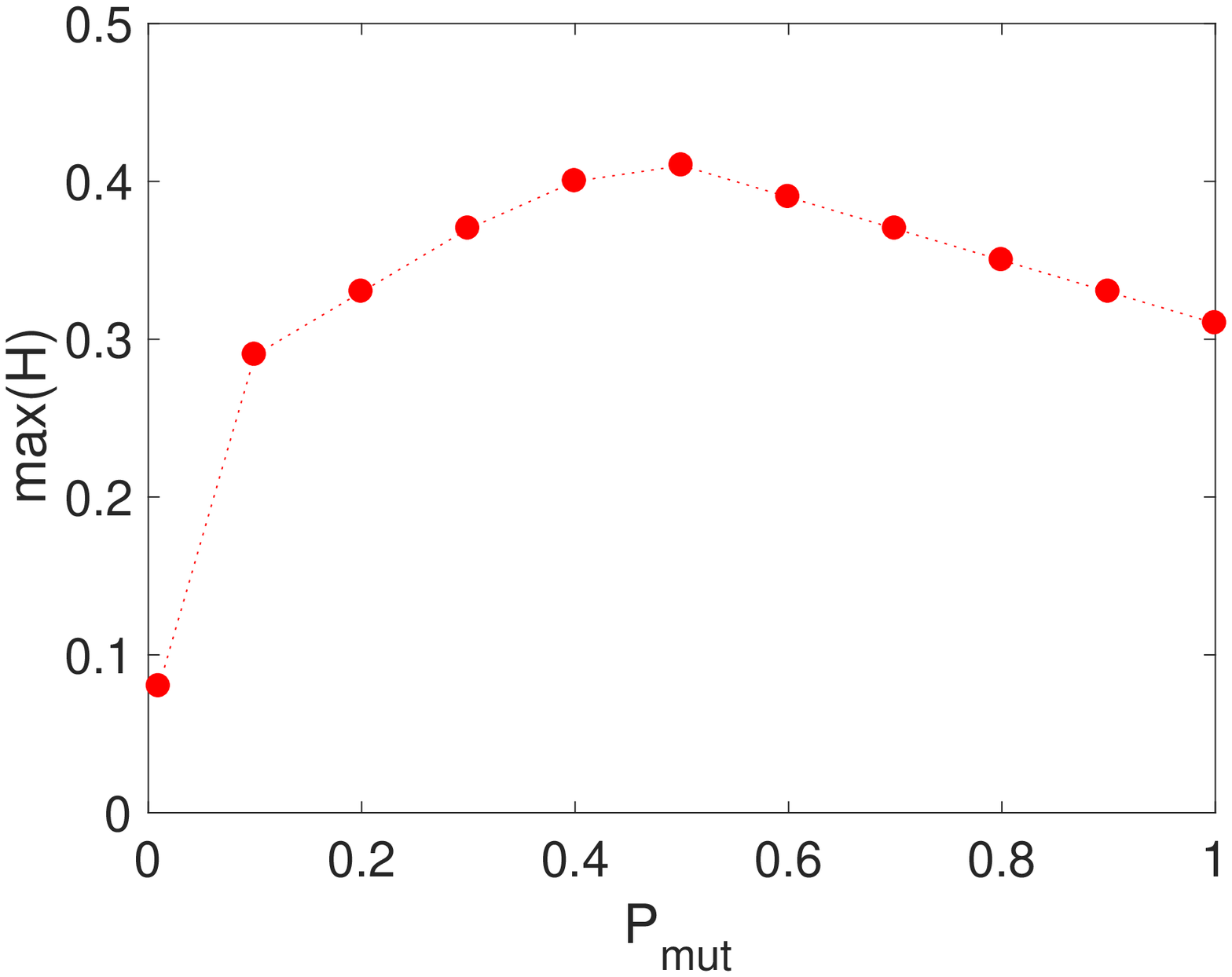}
	\caption{$\max(H(\mathrm{P_{mut}}))$}
\end{subfigure}
\begin{subfigure}[]{0.55\textwidth}
\includegraphics[width=1\textwidth]{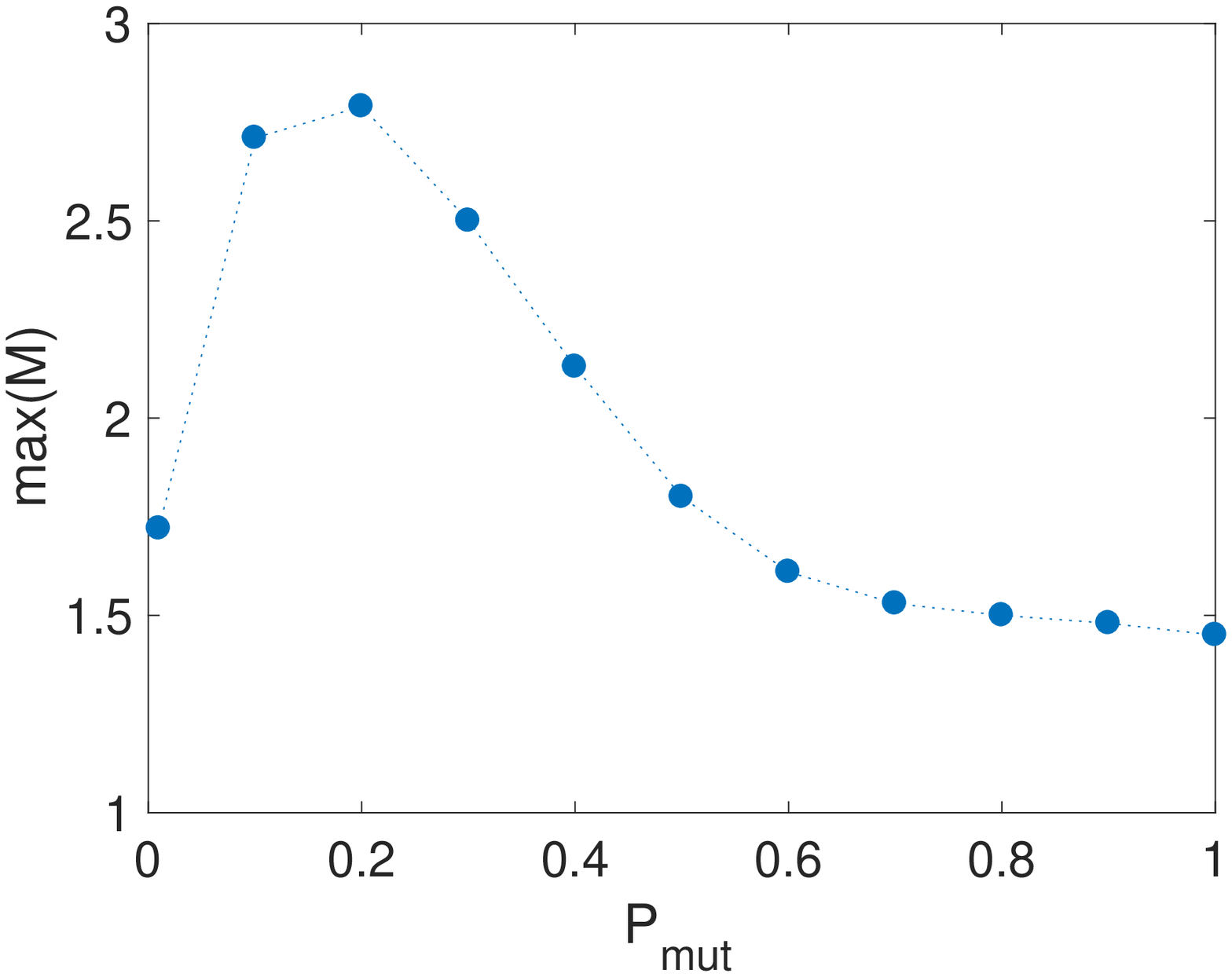}
	\caption{$\max(M(\mathrm{P_{mut}}))$}
\end{subfigure}}
\caption[Dependence of eradication of the secondary lesion from its mutation rate.]{{\bf Dependence of eradication of the secondary lesion from its mutation rate.}
Results for secondary lesions with different mutation rates, when a
combination therapy is administered to the first lesion. 
Panel (A): Total cancer population as function of time for different
$P_{\mathrm{mut}}$. At day $0$, the secondary tumour mass is
phenotypically identical
to the primary lesion. Panel (B) and (C): Maximum values of the Shannon index $H$ and of
the roughness $M$ as functions of the mutation rate. {\it Points in the
Figure represent results from single simulations. Maxima in panels (B)
and (C) have been taken over time for each rate of mutation for the
clones.}} 
\label{au:fig6}
\end{figure}

To investigate how the rate
of mutation of cancer cells affects outcomes on secondary lesions,
we generate a tumour of $5$ x $10^4$ cells with a full mutation capability and apply a combination therapy
(RT and immune boosting) as per the previously introduced
protocols. In other words, the complete secondary lesion before the
start of treatments is composed by the \textit{original} clone.  
Typical results are reported in Fig~\ref{au:fig6}(a) and show that outcomes
do not linearly depend on the rate of mutation. Interestingly, tumour reductions at
day $30$ are larger when the mutation rate is lower, but become
negligible when the rate of mutation is approximately larger than
$P_{\rm mut} = 0.6$, with no relevant change in the overall outcome
for higher rates (see purple, yellow and cyan lines). Also, for rates lower than $P_{\rm mut} = 0.1$ (see
blue and black lines),
different dynamics of eradication can be present, with tumours having
different cell counts after the treatment is administered, although the final result
at day $30$ appears almost identical. 

Heterogeneity tends to increase with the mutation rate, but
its maximum value during therapy is not directly proportional to how
fast the system can mutate. After a given rate of mutation, which is
approximately $P_{\rm mut} \approx 0.5-0.6$, 
mutations do not increase the heterogeneity of the mass. This is
because the fittest
clone usually becomes dominant, its cells outnumber the other
phenotypes and heterogeneity reduces. Fig.~\ref{au:fig6}(b) shows that 
the maximum values of $H$ are a concave function of mutation
probability. Also, the overall 
distribution of phenotypes in the secondary mass is, as discussed
previously for the case of a single cancer, strongly affected by the
morphology.  

Different rates induce different roughness as shown in
Fig.~\ref{au:fig6}(c), where the maximum value of
$M$ reached by the tumour mass is plotted as a function of
mutation rate. This value increases for small rates and reaches a
maximum for a rate $P_{\rm mut} = 0.2$, followed by a sharp
decrease. For rates larger than $P_{\rm mut} = 0.5$, the value does
not change significantly. Tumour morphology of the secondary lesion
influences how the immune system progresses in its attack. In our
simulations, we observe different dynamics of attack carried out by
the immune system, with the tumour being eroded in different ways and
often not in a homogeneous fashion. Nonetheless, a rougher tumour always appears more vulnerable to
immune system attacks because of the degree of infiltration by CTLs it
allows. Even in ``abscopal'' positive outcomes, infiltration plays a major role in the dynamics of erosion and
high $M$ correlates with better results. At the highest mutation rates, roughness is low
because the immune system is not able to recognise
phenotypes that are different from those of the primary mass and kill
them. This results in a fast, unbounded growth of one or two phenotypes that
increase the sphericity of the mass and quickly lower the roughness
value approximately to unity, which corresponds to a spherical
object. This is reflected in the plateau observed in Fig.~\ref{au:fig6} for $P_{\rm mut} > 0.5$.
 
\section{Discussion}\label{sec:discuss}
Cancer and immune cells are complex systems with different
characteristics that also depend on internal and external evolutionary
pressures. In the last decades, improvements on the general
knowledge of these processes have stimulated new
therapeutic strategies which take into account to patients' 
particularities to some degree. The detailed model of
immune interaction described here focuses on the
salient traits of the dynamics and is able to reproduce the
major features of a number of therapeutic interventions.  

An analysis of the effect of different drugs on three prototypical
secondary masses arising from a metastatic breast cancer (not modeled) has been proposed, showing a faithful representation of
 some of the main mechanisms of tumour-immune interactions present in
literature~\cite{piretto2018combination, de2006mixed}. In particular,
we note a significant dependence of the outcomes on the
heterogeneity of the tumour, with higher heterogeneity generally
correlated with negative outcomes confirming biological evidence suggested in Ref.~\cite{burrell2013causes}. 
Indeed, therapies targeting heterogeneous cancer micro environments often show large
rebounds of more resistant tumour cells, that are able to counteract
the action of drugs or boosting in a consistent way. In particular,
one original result from our modeling is that
chemotherapy appears more efficient in a less phenotypically differentiated secondary lesion 
independently of the rate of growth or the dimension of the mass. 
For the reproductive rates considered in this study, 
immune boosting alone is not sufficient to produce full eradication, but rather can trigger
the spread of the more aggressive cells in the body by making the
existing mass more sparse. A well-timed intervention with a combination of
boosting and chemotherapy seems to be the safest of the protocols, 
allowing for a relevant reduction of the mass and preventing the
unbounded growth of the most proliferating cells. Nevertheless, timing of intervention on the secondary lesion can be critical.

One further result of this work is to uncover the importance
of tumour morphology in evolution and fate of secondary lesions. 
The shape of secondary masses conveys important information 
that could be an indicator of successful eradication. For instance, during or
immediately after administration of chemotherapy, our modeling shows
that a high infiltrated tumour is associated with the best outcome. At the same time, a
harsh environment or a high selective pressure tend to generate a tumour that has a
greater roughness and the tendency to spread, as previously noted
~\cite{anderson2006tumor}. This often occurs and persists for many
days after administration of
chemotherapy. In particular, findings for case A suggest that clonal composition of
surviving cells that originated from the beginning tumour colony, and
are later influenced 
by the selected therapies, strongly affect the final outcomes of the
metastasis. Similar dynamics is reported to be present in some types
of tumours that are known to be particularly resistant to therapies~\cite{Kusoglu201980}.

Here we consider one cycle of treatment, but some extra {\it care} should be
taken when multiple cycles are considered, since the immune system
could be further weakened and respond less efficiently.
On the other hand, there is evidence in the medical literature that a
combination of radiotherapy and immunotherapy can provide a positive
effect not only in the area affected by the radiation but also in
other areas of the body. This seems to be due to the release of autologous neoantigens to the
immune system \cite{demaria2015role}, with the overall result of what
it appears to be an individualised tumour vaccine. Our model
captures the effects of this cascading action on a secondary
metastatic mass and confirms that the immune system can act
as a mediator for secondary attacks. In particular, there is evidence
in our findings 
that immune-suppressed hosts or secondary lesions
antigenically unrelated to the treatment area do not show any abscopal 
effect, as experimentally noted. 

There is a growing discussion in the community about the causes of
this rare, positive occurrence on secondary masses. Currently, it
seems that this is the result of a fragile balance between positive
and negative signals activated with the radiation, dependent on the
pre-existing environment and the immunogenicity of the
tumour~\cite{demaria2015role}. These biological elements, alongside
a critical dependence on the dose and the interval between
radiation fractions, contribute to the low occurrence of this effect.
Furthermore and interestingly, a
dependence of the effect on the mutation rate of the cancer clones in
the secondary site is apparent, suggesting that the role of genetic instability at that
site should be investigated more.

Overall, our findings emphasises that the morphology of
the secondary lesions before, during and after the treatment, bears some indications of the rate of success for the treatment.
For the lesions under the detectability threshold, this work suggests that heterogeneity and roughness
are both important quantities. Negative prognosis is linked to the selection of a
poorly immunogenic clone and has been shown to lead to a large,
unbounded regrowth of the tumour. 
It is thus vital to design a protocol that can minimize the immunoediting ability
of cells surviving from therapies or improve the immune system ability to recognise 
and attacks such clones especially when we cannot detect individual,
isolated lesions but only total tumour burden. 
In regards to the latter, the activation of an
``abscopal-like'' response seems to be a strategy in re-calibrating the immune reaction to such
cells. Results are still in their infancy and it is unknown whether such a response can be
elicited and how. One of the ideas
we suggest is to carefully consider the best tumour target to be irradiated:
when possible, it could be advantageous to target a secondary
metastasis antigenically related to the more common lesions in the
body, with a low grade of hypoxia and with a good grade of
immunogenicity. 

\section{Conclusions}
\label{sec:conclusion}
Cellular mutation constitutes one of the causes of negative
outcomes in therapeutic strategies against cancer. Morphology, growth rates and the clonal composition of a tumour mass can, to some extent, be used as predictors of tumour
resistance to a range of anti-cancer therapies and be analysed to combine
treatments to maximize their impacts. In most of the commonly used protocols, the action of the immune system is
crucial. Modern therapies elicit and enhance patient's immune
response, also because of its ability to adapt to, change and modify the
tumour microenvironment. 

The mathematical model we have presented tries to capture the
complexity of tumour-immune dynamics and discuss how
therapies with different scopes, doses and protocols can influence
prognosis for small, solid, secondary tumoural lesions. {\it Even if
  these lesions can be small and not yet detectable, their role can be
  catastrophic for the patient if they are untreated or, as in some
  cases we have shown, the effect of therapies on the primary tumour
  can lead to the selection of more aggressive and resilient clones in
  the secondary lesions.} Ideally,
individual, evidence-based modelling
might provide a fast, reliable and patient-centred way to test and find optimum control of 
protocols {\it in vivo}.  

Our findings suggest that the success of synergistic therapies is
strongly influenced by the phenotypical composition of all the lesions,
alongside their mutation rates and immunogenic properties. 
Effective strategies that can ``normalise'' the
microenvironment \cite{jain2005normalization} and will try to 
limit tumour clonal mutation could be trialled to improve prognosis.  
Therapies that target slowly proliferating and undifferentiated cells
can also become viable in the future.

\section{Acknowledgements}
EP, PSK and FF gratefully acknowledge support for this work through
the Australian Research Council Discovery Project DP180101512,
“Dynamical systems theory and mathematical modelling of viral
infections”. 

\section{Code and Data availability}
The code used and the produced data  for the simulations discussed in this work are available upon request.






\bibliographystyle{model1-num-names}
\bibliography{references.bib}










\end{document}